\documentclass{amsart}
\usepackage{amssymb,amscd}
\DeclareMathSymbol{\twoheadrightarrow} {\mathrel}{AMSa}{"10}

\def\Q{{\mathbf Q}}
\def\Z{{\mathbf Z}}
\def\C{{\mathbf C}}
\def\R{{\mathbf R}}

\def\i{{\mathbf i}}
\def\j{{\mathbf j}}
\def\k{{\mathbf k}}

\def\CC{{\mathfrak C}}

\def\A{{\mathbf A}}

\def\Gal{\mathrm{Gal}}

\def\Is{\mathrm{Isog}}
\def\tr{\mathrm{tr}}

\def\Hdg{\mathrm{Hdg}}
      \def\SU{\mathrm{SU}}
      \def\UU{\mathrm{U}}
\def\End{\mathrm{End}}
\def\Aut{\mathrm{Aut}}
\def\Hom{\mathrm{Hom}}
\def\Mor{\mathrm{Mor}}

\def\I{\mathrm{Id}}
\def\J{{\mathcal J}}

\def\fchar{\mathrm{char}}

        \def\discr{\mathrm{discr}}
\def\GL{\mathrm{GL}}

\def\ZZ{\mathcal{Z}}

        \def\Is{\mathrm{Isog}}

\def\M{\mathrm{M}}
\def\dim{\mathrm{dim}}

\def\P{{\mathbf P}}

\def\X{{\mathcal X}}

\def\Y{{\mathcal Y}}
\def\U{{\mathcal U}}

\newtheorem{thm}{Theorem}[section]
\newtheorem{lem}[thm]{Lemma}
\newtheorem{cor}[thm]{Corollary}
\newtheorem{prop}[thm]{Proposition}

\theoremstyle{definition}
\newtheorem{defn}[thm]{Definition}
\newtheorem{ex}[thm]{Example}

\newtheorem{sect}[thm]{}

           \newtheorem{rem}[thm]{Remark}

\title[Isogeny classes]
{Isogeny classes of abelian varieties over function fields}
\author[Yuri G. Zarhin]{Yuri G. Zarhin}
\address{Department of Mathematics, Pennsylvania State University,
University Park, PA 16802, USA}

 \email{zarhin\char`\@math.psu.edu}
\begin{document}
\begin{abstract}
We study finiteness problems for isogeny classes of abelian varieties over an
algebraic function field $K$ in one variable over the field of complex numbers.
In particular, we construct explicitly a non-isotrivial  absolutely simple
abelian fourfold $X$ over a certain $K$ such that the isogeny class of $X
\times X$ contains infinitely many mutually non-isomorphic principally
polarized abelian varieties. (Such examples do not exist when the ground field
is finitely generated over its prime subfield.)

2000 Mathematics Subject Classification  14K02 (primary), 14D07 (secondary)
\end{abstract}

\maketitle
\section{Introduction}
Let $K$ be a field, $\bar{K}$ its separable closure,
$\Gal(K)=\Gal(\bar{K}/K)$ the (absolute) Galois group of $K$. Let
$X$ be an abelian variety over $K$, $\End_K(X)$ its ring of
$K$-endomorphisms and $\End_K^0(X):=\End_K(X)\otimes\Q$.

Let us consider the set $\Is(X,K)$ of $K$-isomorphism classes of
abelian varieties $Y$ over $K$ such that there exists a
$K$-isogeny $Y \to X$.
If $\ell$ is a prime different from $\fchar(K)$ then let us
consider the set $\Is(X,K,\ell)$ of $K$-isomorphism classes of
abelian varieties $Y$ over $K$ such that there exists an
$\ell$-isogeny $Y \to X$ that is defined over $K$. (Recall that
$\ell$-isogeny is an isogeny, whose degree is a power of $\ell$.)
If $m$ is a
positive integer then we write $\Is_m(X,K,\ell)$ for the subset of
$\Is(X,K,\ell)$ that consists of all (isomorphism) classes of $Y$
with a $K$-polarization of degree $m$.

When $K$ is finitely generated over its prime subfield, Tate
\cite{Tate} conjectured  the finiteness of  $\Is_m(X,K,\ell)$.
This conjecture played a crucial role in the proof of Tate's
conjecture on homomorphisms and the semisimplicity of the Galois
representations in the Tate modules of abelian varieties
\cite{Tate,Serre,Parshin,ZarhinFA,ZarhinIz,ZarhinMZ1,Faltings1,Faltings2,MB,ZarhinP}.

If $K$ is finitely generated over the field $\Q$ of rational
numbers then it was proven by Faltings \cite{Faltings2} that
$\Is(X,K)$ is finite. (See also \cite{Faltings1,ZarhinIn}.) When
$K$ is finite, the finiteness of $\Is(X,K)$ is proven in
\cite{ZarhinMZ2}. (See
\cite{Parshin,ZarhinFA,ZarhinMZ1,ZarhinMZ2,MB} for a discussion of
the case when $K$ is infinite but finitely generated over a finite
field.)

The aim of this note is to discuss the situation when $K$ is an algebraic
function field in one variable over $\C$. (In other words, $K$ is finitely
generated
 and has transcendence degree $1$ over $\C$.) In order to state our results, notice
 that one may view $K$ as
the field of rational functions on a suitable irreducible smooth algebraic
complex curve $S$ such that $X$ becomes the generic fiber of an abelian scheme
$f:\X\to S$. Let $s\in S$ be a complex point of $S$, let $\X_s$ be the fiber of
$f$ at $s$ and $H_1(\X_s,\Q)$ its
 first rational homology group. (Recall that $\X_s$ is a complex abelian variety, whose dimension
  coincides with $\dim(X)$.) There is the natural global monodromy representation
  \[\pi_1(S,s) \to \Aut(H_1(\X_s,\Q))\]
  of the fundamental group $\pi_1(S,s)$ of $S$. Deligne \cite{Deligne} proved that
  this representation is completely reducible and therefore its centralizer
  $D_f:=\End_{\pi_1(S,s)}(H_1(\X_s,\Q))$ is a (finite-dimensional) semisimple $\Q$-algebra. On
  the other hand, every $u\in \End_K(X)$ extends to an endomorphism of $\X$ and therefore induces
  a certain endomorphism $u_s$ of $\X_s$. This gives rise to the embeddings
  \[\End_K(X) \hookrightarrow \End(\X_s) \hookrightarrow
  \End_{\Q}(H_1(\X_s,\Q)),\]
  whose composition extends by $\Q$-linearity to the embedding
  \[\End_K^0(X)\hookrightarrow \End_{\Q}(H_1(\X_s,\Q)),\]
  whose image lies in $D_f$. Further, we identify $\End_K^0(X)$ with its image in $D_f$.
  Our main result (Theorem \ref{main0}) may be restated as follows.

  {\sl The set $\Is(X,K)$ is infinite if and only if $D_f\ne \End_K^0(X)$}.

  If $D_f\ne \End_K^0(X)$ then the set $\Is(X,K,\ell)$ is infinite for all but finitely many
   primes $\ell$ (see Theorem \ref{main}). In addition,
   if $X$ is principally polarized over $K$ then the set
 $\Is_1(X^2,K,\ell)$ is infinite for all but finitely many primes $\ell$ congruent to $1$ modulo $4$
 (see Corollary \ref{principal}).

 In order to describe other results of this paper (Sect. \ref{dim4}), let us further assume that
 every homomorphism over $\overline{\C(S)}$ between $X$ and
``constant" abelian varieties (defined over $\C$) is zero. If $\dim(X)\le 3$
 then Deligne \cite{Deligne} proved that  $D_f= \End_K^0(X)$: this implies that
in this case  $\Is(X,K)$ is finite. On the other hand, Faltings \cite{Faltings0}
constructed a four-dimensional $X$ with $D_f \ne \End_K^0(X)$; in his example(s)
$X$ is the generic fiber of an universal family of abelian fourfolds
with level $n\ge 3$ structure over a Shimura curve.

We prove that if $\dim(X)=4$, all endomorphisms of $X$ are defined over $K$ and $\Is(X,K)$ is infinite
 then $\End_{K}^0(X)$ is a CM-field of degree $4$ (see Theorem \ref{dim4s}).
Almost conversely,  if $X$ is a fourfold, $\End_K^0(X)$ is a CM-field of degree $4$
and all  endomorphisms of $X$ are defined over $K$ then there exists
 a finite algebraic extension $L/K$ such
that $\Is(X\times_{K}L, L)$ is infinite (see Theorem \ref{foursuf}).
A rather explicit example of the latter case (with infinite $\Is(X\times_{K}L, L)$)
is provided by the field
$K=\C(\lambda)$ of rational functions in independent variable $\lambda$,
  the jacobian $X$ of the genus $4$ curve $y^5=x(x-1)(x-\lambda)$ and
the overfield $L=\C(\sqrt[10]{\lambda},\sqrt[10]{\lambda-1})$ (see Example \ref{mainex}).

The paper is organized as follows. Section \ref{prelim} contains basic notation
and useful facts about abelian varieties.
In Section \ref{monod} we
discuss abelian schemes over curves and corresponding monodromy
representations.
 In Section \ref{mr} we state the main results. In Sections
\ref{AS} and \ref{isab} we discuss non-isotrivial abelian schemes
and isogenies of abelian schemes respectively. The next two
sections contains the proofs of main results. Section \ref{append}
contains auxiliary results about quaternions. In Section
\ref{dim4} we deal with isogeny classes of four-dimensional
abelian varieties.

{\bf Acknowledgements}. This paper is based on beautiful ideas and
results of J. Tate \cite{Tate} and P. Deligne \cite{Deligne}. I
learned about the Tate conjecture in the early 1970s from a paper
of A. N. Parshin \cite{Parshin}. I am
 grateful to Minhyong Kim  for  stimulating questions and to B. van
 Geemen  for useful comments. The final version of this paper was written during
 my stay at the Steklov Mathematical Institute of the Russian Academy of
 Sciences (Moscow); I am grateful to its Department of Algebra for the
 hospitality.

 \section{Abelian varieties}
 \label{prelim}
Let $X$ be an abelian variety of positive dimension over a field $K$
 of characteristic zero.
If $n$ is a
positive integer
 then we
write $X_n$ for the kernel of multiplication by $n$ in $X(K_s)$. It is
well-known \cite{Mumford} that $X_n$ ia a free $\Z/n\Z$-module of rank
$2\dim(X)$; it is also a Galois submodule in $X(\bar{K})$. We write $K(X_n)$
for the field of definition of all points of order $n$; clearly, $K(X_n)/K$ is
a finite Galois extension, whose Galois group is canonically identified with
the image of $\Gal(K)$ in $\Aut(X_n)$. We write $\I_X$ for the identity
automorphism of $X$. We write $\End(X)$ for the ring of all
$\bar{K}$-endomorphisms of $X$ and $\End^0(X)$ for the corresponding
$\Q$-algebra $\End(X)\otimes\Q$. We have \[\Z\cdot \I_X \subset \End(X)\subset
\End^0(X).\]
\begin{rem}
\label{isogeny} Let $X'$ be an abelian variety over $K$. If $X$
and  $X'$  are $K$-isogenous then
every $K$-isogeny $X \to X'$ gives rise to
 a bijection between
$\Is(X,K$ and $\Is(X,K')$. In
particular, the sets $\Is(X,K)$ and $\Is(X',K)$ are either both
finite or both infinite.
\end{rem}

Let $\ell$ be a prime.

\begin{lem}
\label{onemod4}
Suppose that $\ell$ is congruent to $1$ modulo $4$. If $X$ admits a principal
polarization over $K$ and
$\Is(X,K,\ell)$ is infinite
then $\Is_1(X^2,K,\ell)$ is also infinite.
\end{lem}

\begin{proof}
Clearly, every abelian
 variety $Y$ over $K$ that admits an $\ell$-isogeny $Y\to X$
 over $\C(S)$ also admits a $K$-polarization, whose degree is
 a power of $\ell$. In addition, there exists an $\ell$-isogeny $X\to Y$
 over $K$;
 its {\sl dual}  $Y^t\to X^t\cong X$ is an $\ell$-isogeny that is also defined over
 $K$.
 Since $\ell$ is congruent to $1$ modulo $4$ then
 $\sqrt{-1}\in \Z_{\ell}$ and it follows from Remark 5.3.1 on pp.
 314--315 of \cite{ZarhinIn} that $Y \times Y^{t}$ is principally
 polarized over $K$. Clearly, there is an $\ell$-isogeny $Y \times Y^{t}\to X\times X=X^2$
 that is defined over $K$. On the other hand, for any given abelian variety
 $Z$ over $K$
 there are, up to an isomorphism, only finitely many abelian
 varieties over $K$
  that are isomorphic over $K$ to an abelian subvariety of $Z$  \cite{LOZ}. It follows that
 $\Is_1(X^2,K,\ell)$ is {\sl infinite}.
\end{proof}

\

\begin{rem}
\label{twist}
 Let $L/K$ be a finite Galois extension with Galois group
 $\Gal(L/K)$. Let $Y$ be an abelian variety over $K$ and
 $\Aut_L(Y)$ be its group of $L$-automorphisms. It is well-known
 that $\Aut_L(Y)$ is an {\sl arithmetic} group; it follows from a
 theorem of Borel-Serre  \cite{BS} that the corresponding first noncommutative Galois cohomology set
 $H^1(\Gal(L/K), \Aut_L(Y))$ is finite, i.e., the set of
 ($K$-isomorphism classes of)
 $L/K$-forms of $Y$ is finite. This implies that if $\Is(X\times_{K}L,L)$ (resp.
 $\Is(X\times_{K}L,L,\ell)$)
 is finite then  $\Is(X,K)$ (resp. $\Is(X,K,\ell)$) is also finite \cite{ZarhinP}.
\end{rem}

Abusing notation, we sometime write $\Is(X,L,\ell)$, $\Is(X,L,\ell)$ and
$\Is_m(X,L,\ell)$) instead of $\Is(X\times_{K}L,L,\ell)$,
$\Is(X\times_{K}L,L,\ell)$ and $\Is_m(X\times_{K}L,L,\ell)$ respectively.

\section{Abelian schemes and monodromy representations}
\label{monod}

 Let $S$ be an irreducible smooth (but not necessarily projective) algebraic curve over
$\C$. We write $\C(S)$ for the field of rational functions on $S$
and $\overline{\C(S)}$ for its algebraic closure. Let $f:\X \to S$
be a polarized abelian scheme of positive relative dimension $d$
over $ S$. Let $\eta$ be the generic point of $S$ and  $X$ be the
generic fiber; it is a $d$-dimensional abelian variety over
$k(\eta)=\C(S)$.

An abelian variety $Z$ over $\C(S)$ is called {\sl constant} if
there exists an abelian variety $W$ over $\C$ such that $Z$ is
isomorphic to $W\times_{\C}\C(S)$ over $\C(S)$.

Recall that $\X$ is called {\sl isotrivial} if there exists a
finite  \'etale cover $S'\to S$  (with non-empty $S'$) such that
the pullback $\X_{S'}=\X\times_S S'$ is a constant abelian scheme
over $S'$. We say that $\X$ is {\sl weakly isotrivial} if (under
the same assumptions on $S'\to S$) $\X_{S'}$ contains a non-zero
constant abelian subscheme.

\begin{thm}
\label{iso}
\begin{itemize}
 \item[(i)] $\X$ is isotrivial if and only if there
exists an abelian variety $W$ over $\C$ such that the
$\C(S)$-abelian varieties $X$ and $W\times_{\C}\C(S)$ are
isomorphic over $\overline{\C(S)}$. \item[(ii)]$\X$ is weakly
isotrivial if and only if there exist an abelian variety $W$ of
positive dimension over $\C$  and an abelian subvariety $Z\subset
X\times_{\C(S)}\overline{\C(S)}$ such that the
$\overline{\C(S)}$-abelian varieties
$W\times_{\C}\overline{\C(S)}$ and $Z$ are isomorphic.
\end{itemize}
\end{thm}

We prove Theorem \ref{iso} in Section \ref{AS}.

We write $R_1 f_*\Z$ for the corresponding local system of the first integral
homology groups $H_1(\X_s,\Z)$ of the fibres $\X_s=f^{-1}(s)$ on $S$. If $g: \Y
\to S$ is another polarized abelian scheme over $S$ with generic fiber $f$ then
the natural map \[\Hom_S(\Y,\X) \to \Hom_{\C(S)}(Y,X)\] is bijective. We write
$R_1 f_*\Q$ for the corresponding local system  $R_1 f_*\Z\otimes\Q$ of the
first rational homology groups $H_1(\X_s,\Q)=H_1(\X_s,\Z)\otimes\Q$ of the
fibres $\X_s=f^{-1}(s)$ on $S$ and $\pi_1(S,s) \to \Aut(H_1(\X_s,\Z))\subset
\Aut(H_1(\X_s,\Q))$ for the corresponding monodromy representation.

\begin{rem}
\label{kz} An abelian scheme $f:\X\to S$ is isotrivial if and only if the image
of the monodromy representation $\pi_1(S,s)\to\Aut(H_1(\X_s,\Z))$ is finite
\cite{Katz} (see also \cite[Sect. 4.1.3.3]{Deligne}).
\end{rem}

\begin{sect} {\bf Rigidity and specialization}.
\label{rigid} It follows from the rigidity lemma \cite[Ch. 6, Sect. 1, Cor.
6.2]{GIT} that the natural homomorphism \[\Hom_S(\Y,\X)\to\Hom(\Y_s,\X_s),
\quad u \mapsto u_s\] is injective. In particular, the natural ring
homomorphism $\End_S(\X) \to \End(\X_s)$ is also injective. It is well-known
that the natural homomorphisms
\[\Hom(\Y_s,\X_s)\to \Hom(H_1(\Y_s,\Z),H_1(\X_s,\Z)), \ \End(\X_s)\to
\End(H_1(\X_s,\Z))\] are embeddings. Taking the compositions, we get the
embeddings \[\Hom_S(\Y,\X)\to \Hom(H_1(\Y_s,\Z),H_1(\X_s,\Z)), \End_S(\X) \to
\End(H_1(\X_s,\Z));\] the images of  $\Hom_S(\Y,\X)$ and of  $\End_S(\X)$ lie
in $\Hom_{\pi_1(S,s)}(H_1(\Y_s,\Z),H_1(\X_s,\Z))$ and
$\End_{\pi_1(S,s)}(H_1(\X_s,\Z))$ respectively. Further we will identify
$\End_S(\X)$ with its image in $\End_{\pi_1(S,s)}(H_1(\X_s,\Z))\subset
\End(H_1(\X_s,\Z))$ and $\Hom_S(\Y,\X)$ with its image in
$\Hom_{\pi_1(S,s)}(H_1(\Y_s,\Z),H_1(\X_s,\Z))$
 respectively.

 We have
\[\Hom^0_S(\Y,\X):=\Hom_S(\Y,\X)\otimes\Q\subset\]
\[\Hom_{\pi_1(S,s)}(H_1(\Y_s,\Z),H_1(\X_s,\Z))\otimes\Q=
\Hom_{\pi_1(S,s)}(H_1(\Y_s,\Q),H_1(\X_s,\Q)),\]
\[\End^0_S(\X):=\End_S(\X)\otimes\Q\subset
\End_{\pi_1(S,s)}(H_1(\X_s,\Z))\otimes\Q= \End_{\pi_1(S,s)}(H_1(\X_s,\Q)).\]
Notice that a theorem of Grothendieck \cite[p. 60]{Grothendieck} (see also
\cite[Sect. 4.1.3.2] {Deligne}) implies that
\[\Hom_S(\Y,\X)=\Hom(\Y_s,\X_s)\bigcap
\Hom_{\pi_1(S,s)}(H_1(\Y_s,\Z),H_1(\X_s,\Z)),\]
\[\End_S(\X)=\End(\X_s)\bigcap\End_{\pi_1(S,s)}(H_1(\X_s,\Z)).\] It follows that
\[\Hom^0_S(\Y,\X)=\Hom^0(\Y_s,\X_s)\bigcap
\Hom_{\pi_1(S,s)}(H_1(\Y_s,\Q),H_1(\X_s,\Q)).\] Here
$\Hom^0(\Y_s,\X_s)=\Hom(\Y_s,\X_s)\otimes\Q \subset
\Hom_{\Q}(H_1(\Y_s,\Q),H_1(\X_s,\Q))$. In particular,
\[\End_S^0(X)=\End^0(\X_s)\bigcap\End_{\pi_1(S,s)}(H_1(\X_s,\Q)).\]
\end{sect}

\begin{sect}
\label{deligneT}
 Let $\Gamma_s$ be the image of $\pi_1(S,s)$ in
$\Aut(H_1(\X_s,\Q))$, let  $G_s$ is the Zariski closure of
$\Gamma_s$ in $\GL_{\Q}(H_1(\X_s,\Q))$ and $G_s^0$ its identity
component. By a theorem of Deligne \cite[Cor. 4.2.9]{Deligne},
$G_s^0$ is a {\sl semisimple} algebraic $\Q$-group.
      Recall \cite[Ch. 1, Sect. 1.2]{BorelLin} that  $G_s^0$ has
              finite index in $G_s$.
It
follows that the intersection $\Gamma_s^0:=\Gamma_s\bigcap G_s^0$ is a normal
subgroup of finite index in $\Gamma_s$ and this index divides $[G_s:G_s^0]$.

The groups $\Gamma_s$ and $G_s$ do not depend on $s$ in the
 following sense. Recall that $S$ is arcwise connected with respect to the complex topology.
 If $t$ is another point of $S$ then every path $\gamma$ in $S$ from $s$ to
 $t$ defines an isomorphism $\gamma_{*}:H_1(\X_s,\Z)\cong
 H_1(\X_t,\Z)$ such that
 $\Gamma_t=\gamma_{*} \Gamma_s \gamma_{*}^{-1}$ and therefore
 \[G_t=\gamma_{*}G_s \gamma_{*}^{-1}, \ G_t^0=\gamma_{*}G_s^0
 \gamma_{*}^{-1}.\]
\end{sect}

  Let $\Hdg(\X_s) \subset
\GL_{\Q}(H_1(\X_s,\Q))$ be the {\sl Hodge group} of  $\X_s$ \cite{Mumford} (see
also \cite{ZarhinIzv,MZ,MZ2}). Recall that $\Hdg(\X_s)$  is a connected
reductive algebraic $\Q$-group and its centralizer
$\End_{\Hdg(\X_s)}(H_1(\X_s,\Q))$ in $\End_{\Q}(H_1(\X_s,\Q))$ coincides with
$\End^0(\X_s)$.

The following assertion is a special case of a theorem of Deligne
\cite[Th. 7.3]{ZarhinIzv} (see also \cite[Prop. 7.5]{DeligneK3},
\cite{Masser}).

\begin{thm}
\label{deligneH} For all $s \in S$ outside a countable set, $G_s^0$ is a normal
algebraic subgroup in $\Hdg(\X_s)$.
\end{thm}

\begin{defn}
We say that a point $s \in S$ is {\sl in general position} (with respect to $\X
\to S$ if $G_s^0$ is a normal algebraic subgroup in $\Hdg(\X_s)$. It follows
from Theorem \ref{deligneH} that every $s$ outside a countable set is in
general position. In particular, a point in general position always does exist.
\end{defn}

\begin{cor}
\label{endoS} Suppose that $s$ is in general position and $G_s$ is connected.
Then \[\End^0(\X_s)=\End_S^0(\X)=\End^0_{\C(S)}(X).\] If, in addition, $\X$ is
not isotrivial and $\Hdg(\X_s)$ is a $\Q$-simple algebraic group then
$\End^0_{\C(S)}(X)=\End_S^0(\X)=\End_{\pi_1(S,s)}(H_1(\X_s,\Q))$.
\end{cor}

\begin{proof}
We have $G_s=G_s^0\subset \Hdg(\X_s)$. Clearly,
\[\End_{\pi_1(S,s)}(H_1(\X_s,\Q))=\End_{G_s}(H_1(\X_s,\Q))\supset
\End_{\Hdg(\X_s)}(H_1(\X_s,\Q))=\End^0(\X_s).\] It follows that
$\End_S^0(\X)=\End^0(\X_s)\bigcap\End_{\pi_1(S,s)}(H_1(\X_s,\Q))=\End^0(\X_s)$.
Assume now that $\X$ is not isotrivial and $\Hdg(\X_s)$ is $\Q$-simple. By
Remark \ref{kz},  $\Gamma_s$ is infinite and therefore $G_s$ has positive
dimension. Since $s$ is in general position, $G_s$ is normal in $\Hdg(\X_s)$.
Now the simplicity of $\Hdg(\X_s)$ implies that $\Hdg(\X_s)=G_s$ and therefore
\[\End^0(\X_s)=\End_{G_s}(H_1(\X_s,\Q))=\End_{\pi_1(S,s)}(H_1(\X_s,\Q))=
\End_S^0(X)=\End_{\C(S)}(X).\]
\end{proof}

\begin{sect} {\bf Base change}.
\label{cover} Let $S' \to S$ be a finite \'etale cover and $\X'=X\times_S S'$
the corresponding abelian $S$-scheme. If a point $s'\in S'$ lies above $s \in
S$ then (in obvious notations)
\[\X'_{s'}=\X_s, \
H_1(\X_s,\Z)=H_1(\X'_{s'},\Z),\ H_1(\X_s,\Q)=H_1(\X'_{s'},\Q),\
\Hdg(\X'_{s'})=\Hdg(\X_s).\]
 The fundamental group $\pi_1(S',s)$ of $S'$ is a
subgroup of finite index in $\pi_1(S,s)$ and therefore $\Gamma_{s'}$ is a
subgroup of finite index in $\Gamma_s$. It follows that $G_{s'}$ is a subgroup
of finite index in $G_s$ ; in particular, $G_{s'}^0=G_{s}^0$. This implies that
$s$ is in general position with respect to $\X$ if and only if $s'$ is in
general position with respect to $\X'$. It also follows that if $G_s$ is
connected then $G_{s'}=G_s$.

On the other hand, let $\Gamma^{0}\subset \Gamma_s$ be a subgroup
of finite index (e.g., $\Gamma^{0}=\Gamma_s^{0}$). Clearly, its
Zariski closure lies between $G_{s}^0$ and $G_{s}$; in particular,
if $\Gamma^{0}\subset\Gamma_s^{0}$ then this closure coincides
with $G_{s}^0$. Let $\pi^{0}\subset \pi_1(S,s)$ be the preimage of
$\Gamma^{0}$: it is a  subgroup of finite index in $\pi_1(S,s)$.
Let $S^{0} \to S$ be a finite \'etale map of irreducible smooth
algebraic curves attached to $\pi^{0}$ with
$\pi_1(S^{0})=\pi^{0}$. Clearly, the degree of $S^{0} \to S$
coincides with the index $[\Gamma_s:\Gamma^0]$.

Notice that if $\Gamma^0$ is normal in $\Gamma_s$ then $S^0\to S$
is Galois with Galois group $\Gamma_s/\Gamma^0$.

Fix a point $s_0\in S^0$ that lies above $s$. Then the image of the monodromy
representation $\pi_1(S^0,s_0)=\pi^{0} \to \Aut(H_1(\X_s,\Q))$ attached to the
abelian $S^{0}$-scheme $\X\times_S S^0$ coincides with $\Gamma^{0}$. In other
words, $\Gamma_{s_0}=\Gamma^{0}$ and therefore $G_{s_0}=G_s^0$ is connected if
$\Gamma^{0}\subset\Gamma_s^{0}$.
\end{sect}

\begin{ex}
\label{min}
 Let $n \ge 3$ be an integer and $\Gamma^{n}\subset \Gamma_s$
 the kernel of the reduction map modulo $n$
\[\Gamma_s \subset \Aut(H_1(\X_s,\Z)) \twoheadrightarrow
\Aut(H_1(\X_s,\Z/n\Z)).\] Here $H_1(\X_s,\Z/n\Z)=H_1(\X_s,\Z)\otimes\Z/n\Z$ is
the first integral homology group of $\X_s$ with coefficients in $\Z/n\Z$.
Clearly, $\Gamma^{n}$ is a (normal) subgroup of finite index in $\Gamma_s$. On
the other hand, since $\Gamma^{n}\subset 1+ n\End(H_1(\X_s,\Z))$, Zariski
closure of $\Gamma^{n}$ is connected \cite[Prop. 2.6 ]{SZcomp}. This implies
that this closure coincides with $G_{s}^0$.
\end{ex}

\begin{rem}
\label{endostab}
\begin{itemize}
\item[(i)] It is known \cite{Silverberg} that all endomorphisms of $X$ are
defined over $L:=\C(S)(X_3)$. Clearly, $L/\C(S)$ is a finite Galois extension
and all points of order $3$ on $X$  are $L$-rational. It follows from
N\'eron-Ogg-Shafarevich criterion \cite{SerreTate} that $L/\C(S)$ is unramified
at all points of $S$. This implies that if $S'$ is the normalization of $S$ in
$K$ then $S'\to S$ is finite \'etale. Clearly, $\C(S')=L$. \item[(ii)] Suppose
that $G_s$ is connected. Then $\End^0(X)=\End_{\C(S)}^0(X)$, i.e., all
endomorphisms of $X$ are defined over $\C(S)$. Indeed, first,
\[\End_{\pi_1(S,s)}(H_1(\X_s,\Q))=\End_{\Gamma_s}(H_1(\X_s,\Q))=\End_{G_s}(H_1(\X_s,\Q))\]
and therefore $\End_{\C(S)}^0(X)=\End^0(\X_s)\bigcap \End_{G_s}(H_1(\X_s,\Q))$.
 Second, pick a point $s'\in S'$ that lies above $s$.
We have $\X'_{s'}=\X_s,\ G_s=G_{s'}$ and
\[\End^0(X)=\End_{\C(S')}(X)=\End^0(\X'_{s'})\bigcap\End_{\pi_1(S',s')}(H_1(\X'_{s'},\Q))=\]
\[\End^0(\X'_{s'})\bigcap\End_{G_{s'}}(H_1(\X'_{s'},\Q))= \End^0(\X_s)\bigcap
\End_{G_s}(H_1(\X_s,\Q))=\End_{\C(S)}(X).\]
\end{itemize}
\end{rem}

\begin{sect} {\bf Isogenies}. \label{isogenySS} It is well-known that $u\in
\Hom_{\C(S)}(Y,X)$ is an isogeny (resp. an $\ell$-isogeny) if and
only if there exist $v\in \Hom_{\C(S)}(X,Y)$ and a positive
integer $n$ such that the compositions $uv$ and $vu$ are
multiplications by $n$ (resp. by $\ell^n$) in $X$ and $Y$
respectively.

We say that $u\in\Hom_S(\Y,\X)$ is an isogeny (resp. an
$\ell$-isogeny) of abelian schemes if the induced homomorphism of
generic fibers $u_{\eta}:Y\to X$ is an isogeny of abelian
varieties.

Clearly, $u$ is an isogeny (resp. an $\ell$-isogeny) of abelian
schemes if and only if there exist $v\in \Hom_S(X,Y)$ and a
positive integer $n$ such that the compositions $uv$ and $vu$ are
multiplications by $n$ (resp. by $\ell^n$) in $\X$ and $\Y$
respectively (see \cite[Sect. 7.3]{Neron}).

If $u \in \Hom_S(\Y,\X)$ is an isogeny (resp. an $\ell$-isogeny) of abelian
schemes then it is clear that the induced homomorphism $u_s:\Y_s\to \X_s$ is an
isogeny (resp. an $\ell$-isogeny) of the corresponding complex abelian
varieties \cite[Sect. 7.3]{Neron}.

Conversely, suppose that $u \in \Hom_S(\Y,\X)$ and assume that $u_s \in
\Hom(\Y_s,\X_s)$ is an isogeny (resp. an $\ell$-isogeny). Then there exist an
isogeny $w_s:\X_s\to \Y_s$ (resp. an $\ell$-isogeny)  and a positive integer
$n$ such that the compositions $u_s w_s$ and $w_s u_s$ are multiplications by
$n$ (resp. by $\ell^n$) in $\X_s$ and $\Y_s$ respectively. Since $u_s\in
\Hom_{\pi_1(S,s)}(H_1(\Y_s,\Z),H_1(\X_s,\Z))$ and therefore \[u_s^{-1}\in
\Hom_{\pi_1(S,s)}(H_1(\X_s,\Z)\otimes\Q,H_1(\Y_s,\Z))\otimes\Q),\] it follows
that $w_s\in\Hom_{\pi_1(S,s)}(H_1(\X_s,\Z),H_1(\Y_s,\Z))$. This implies that
\[w_s\in \Hom(\X_s,\Y_s)\bigcap \Hom_{\pi_1(S,s)}(H_1(\X_s,\Z),H_1(\Y_s,\Z))\]
and therefore there exists $v\in \Hom_S(\X,\Y)$ with $v_s=w_s$. It follows that
$(uv)_s$ and $(vu)_s$ are multiplications by $n$ (resp. by $\ell^n$) in $\X_s$
and $\Y_s$ respectively. By the rigidity lemma (Sect. \ref{rigid}), $uv$ and
$vu$ are multiplications by $n$ (resp. by $\ell^n$) in $\X$ and $\Y$
respectively. This implies that $u$ and $v$ are isogenies (resp.
$\ell$-isogenies) of abelian schemes.
\end{sect}

\begin{sect} {\bf Semisimplicity}.
\label{semisimple}
 Recall \cite[Sect. 4.2]{Deligne} that the
{\sl monodromy} representation \[\pi_1(S,s)\to
\Aut(H_1(\X_s,\Z)\subset\Aut(H_1(\X_s,\Z)\otimes\Q))=\Aut_{\Q}(H_1(\X_s,\Q))\]
is completely reducible and  therefore its centralizer
\[D=D_f:=\End_{\pi_1(S,s)}(H_1(\X_s,\Q))\subset\End_{\Q}(H_1(\X_s,\Q))\] is a
finite-dimensional {\sl semisimple} $\Q$-algebra. In addition, the {\sl center}
$E$ of $D$ lies in $\End_S(\X)\otimes\Q$ \cite[Sect. 4.4.7]{Deligne}. It
follows that $E$ lies in the {\sl center} $C$ of $\End^0_S(\X)$. Notice that
$\End_S(\X)\otimes\Q=\End_{\C(S)}(X)\otimes\Q$. So, we have the inclusion of
finite-dimensional semisimple $\Q$-algebras
\[\End_{\C(S)}(X)\otimes\Q=:\End^0_S(\X)\subset D \subset
\End_{\Q}(H_1(\X_s,\Q));\] in addition, $E$ lies in $C$. It follows easily that
if $D$ is commutative then $\End^0_S(\X)= D$.
\end{sect}

\begin{rem}
\label{isotype} Assume that  $X$ is {\sl simple}. Then:
\begin{itemize}
\item[(i)] $\End^0_{\C(S)}(X)$ is a division algebra. Since  the center $E$ of
$D_f$ lies in $\End^0_{\C(S)}(X)$ (Sect. \ref{semisimple}), it has no zero
divisors. This implies that $E$ is a field and therefore $D$ is a central
simple $E$-algebra. It follows that the $H_1(\X_s,\Q)$ is an isotypic
$\pi_1(S,s)$-module, i.e., is either simple or isomorphic to a direct sum of
several copies of a simple module. \item[(ii)] Notice that $\End^0_S(\X)= D$ if
and only if $\End_S(\X)=\End_{\pi_1(S,s)}(H_1(\X_s,\Z))$. It is also known
\cite[Cor. 4.4.13]{Deligne} that $\End^0_S(\X)= D$ when
 $d \le 3$  and $\X$ is not weakly isotrivial.
 \end{itemize}
\end{rem}

 \begin{rem}
\label{fal}
  Faltings \cite[Sect. 5]{Faltings0}
 has constructed a principally polarized abelian scheme $f:\X\to S$ with $d=4$ that is
 {\sl not} weakly isotrivial,  $\Gamma_s=\Gamma^{n}$ for a
 certain integer $n \ge 3$ (in notations of Example \ref{min}) and $\End^0_S(\X)\ne D$.
It follows from arguments in Example \ref{min} that in Faltings'
example $G_s$ is connected.
 \end{rem}

\section{Main results}
\label{mr} Our main result is the following statement.

\begin{thm}
\label{main0} Let $f:\X\to S$ be a polarized abelian scheme of
positive relative dimension $d$. Then the following conditions are
equivalent:

\begin{itemize}
\item [(i)] $D_f \ne \End^0_S(\X)$. \item [(ii)] The set
$\Is(X,\C(S))$ is infinite.
\end{itemize}
\end{thm}

Theorem \ref{main0} is an immediate corollary of the following two
statements.

\begin{thm}
\label{finite} Let $f:\X\to S$ be a polarized abelian scheme of
positive relative dimension $d$. If $D_f = \End^0_S(\X)$ then
$\Is(X,\C(S))$ is finite.
\end{thm}

\begin{thm}
\label{main} Let $f:\X\to S$ be a polarized abelian scheme of
positive relative dimension $d$. Suppose that $D=D_f \ne
\End^0_S(\X)$ (in particular, $D$ is noncommutative). If $\ell$ is
a prime such that $D\otimes_{\Q}\Q_{\ell}$ is isomorphic to a
direct sum of matrix algebras over fields then there exist a
positive integer $d_0<2d$ and a sequence $\{Z(m)\}_{m=1}^{\infty}$
of $\C(S)$-abelian varieties $Z(m)$ that enjoy the following
properties:
\begin{itemize}
\item[(i)] For each positive integer $m$ there exists a
$\C(S)$-isogeny $Z(m)\to X$ of degree $\ell^{(2d-d_0) m}$.
\item[(ii)] Let $Y$ be an abelian variety  over $\C(S)$ and
 let $M_Y$ be the set of positive integers $m$ such that $Z(m)$
 is isomorphic to $Y$ over $\C(S)$. Then $M_Y$ is either empty or
 finite. In other words, if
 $M$ is an arbitrary infinite set of positive integers
then there exists an infinite subset $M_0\subset M$ such that for
$m \in M_0$ all $Z(m)$ are mutually non-isomorphic over $\C(S)$.
\end{itemize}
\end{thm}

\begin{rem}
\label{splitD}
Clearly, $D\otimes_{\Q}\Q_{\ell}$ is isomorphic to a direct sum of
matrix algebras over fields  for all but finitely many primes
$\ell$.
\end{rem}

We prove Theorems \ref{main} and \ref{finite} in Sections
\ref{inf} and \ref{fin} respectively.

\begin{cor}
Suppose that a point  $s\in S$ is in general position and $\Hdg(\X_s)$ is
$\Q$-simple. If $\X$ is not weakly isotrivial then for all finite algebraic
extension $L/\C(S)$ the set $\Is(X,L)$ is finite.
\end{cor}

\begin{proof}
Replacing $S$ by its normalization in $L$ and using  (\ref{cover}), we may
assume without loss of generality that $L=K$. Choosing an integer $n \ge 3$,
replacing $S$ by $S_n$ and applying (\ref{cover}) and Remark \ref{twist}, we
may assume without loss of generality that $G_s$ is connected. It follows from
Corollary \ref{endoS} that
\[\End_{\C(S)}(X)=\End_S^0(\X)=\End_{\pi_1(S,s)}(H_1(\X_s,\Q)).\] Now the result
follows from Theorem \ref{main0}.
\end{proof}

\begin{cor}
\label{principal} Suppose that $D_f \ne \End^0_S(\X)$. If $X$
admits a principal polarization over $\C(S)$ then
$\Is_1(X^2,\C(S),\ell)$ is {\sl infinite} for all but finitely
many primes $\ell$  congruent to $1$ modulo $4$.
\end{cor}

\begin{proof}

Let us pick a prime $\ell$ such that $D_f\otimes_{\Q}\Q_{\ell}$ is
a direct sum of
 matrix algebras over fields. It follows from Theorem \ref{main}
 that $\Is(X,\C(S),\ell)$ is {\sl infinite}. Applying Lemma \ref{onemod4},
 we conclude that $\Is_1(X^2,\C(S),\ell)$ is {\sl infinite} if
  $\ell$  is congruent to $1$ modulo $4$. Now the result follows from Remark \ref{splitD}.
\end{proof}

Till the end of this Section we assume that  $K$ is a field of
algebraic functions in one variable over $\C$, i.e., $K$ is
finitely generated and of degree of transcendency $1$ over $\C$.

\begin{thm}
\label{dim3}  Let $X$ be an abelian variety over $K$ of positive
dimension $d$. Suppose that the $\bar{K}/\C$-trace of $X$ is zero.
If $d\le 3$ then $\Is(X,K)$ is finite.
\end{thm}

\begin{proof}[Proof of Theorem \ref{dim3}]
There exists a smooth irreducible algebraic $\C$-curve $S$ such that $K=\C(S)$
and there exists a polarized abelian scheme $\X \to S$ of relative dimension
$d\le 3$, whose generic fiber coincides with $X$ \cite[Sect.1.4, p.20]{Neron}.
Pick a point $s\in S(\C)$. Recall \cite[Ch. VIII, Sect. 3]{Lang} that the
vanishing of the $\overline{\C(S)}/\C$-trace means that $X$ does not contain
over $\overline{\C(S)}$ a non-zero abelian subvariety that is isomorphic over
$\overline{\C(S)}$  to a (constant) abelian variety of the form
$W\times_{\C}\overline{\C(S)}$ where $W$ is a complex abelian variety. By
Theorem \ref{iso}, this means $X$ is not weakly isotrivial. It follows from
results of Deligne (Remark \ref{isotype}(ii)) that
$\End_S(\X)=\End_{\pi_1(S,s)}(H_1(\X_s,\Z))$ and therefore $D_f:
=\End_{\pi_1(S,s)}(H_1(\X_s,\Q))= \End^0_S(\X)$. Now the assertion follows from
Theorem \ref{main}.
\end{proof}

\begin{cor}
\label{prod3} Let $X$ be an abelian variety over $K$ that is
isogenous over $\bar{K}$ to a product of abelian varieties of
dimension $\le 3$. If $\bar{K}/\C$-trace of $X$ is zero then
$\Is(X,K)$ is finite.
\end{cor}

\begin{proof}
Replacing if necessary, $K$ by its finite Galois extention and
using Remarks \ref{twist} anad \ref{isogeny}, we may assume
without loss of generality that $X$ is isomorphic over $K$ to a
product $Y_1\times \cdots \times Y_r$ of abelian varieties $Y_i$'s
of  dimension $\le 3$ over $K$.

There exists a smooth irreducible algebraic $\C$-curve $S$ such
that $K=\C(S)$ and there exist polarized abelian schemes $\Y_i \to
S$, whose generic fiber coincides with $Y_i$ \cite[Sect.1.4,
p.20]{Neron}. If $f: \X\to S$ is the fiber product of $\Y_i$ then
the generic fiber of the abelian $S$-scheme $\X$ coincides with
$X$.

Pick a point $s\in S(\C)$. Since all $\dim(Y_i)\le 3$, it follows from results
of Deligne \cite[Cor. 4.4.13]{Deligne} applied  to all pairs $(\Y_i,\Y_j)$ that
$\End_S(\X)=\End_{\pi_1(S,s)}(H_1(\X_s,\Z))$ and therefore $D_f:
=\End_{\pi_1(S,s)}(H_1(\X_s,\Q))= \End^0_S(\X)$. Now the assertion follows from
Theorem \ref{main}.
\end{proof}

\begin{cor}
\label{dim4ns} Let  $X$ be a four-dimensional abelian variety over
$K$. Suppose that the $\bar{K}/\C$-trace of $X$ is zero. If $X$ is
not absolutely simple then
 $\Is(X,K)$ is finite.
\end{cor}

\begin{proof}
Replacing if necessary, $K$ by its finite Galois extention and
using Remarks \ref{twist} anad \ref{isogeny}, we may assume
without loss of generality that $X$ is isomorphic over $K$ to a
product $Y\times T$ of abelian varieties $Y$ and $T$ of positive
dimension over $K$. Since $4= \dim(X)=\dim(Y)+\dim(T)$, we
conclude that both $\dim(Y)$ and $\dim(T)$ do not exceed $3$. One
has only to apply Corollary \ref{prod3}.
\end{proof}

\begin{ex}
\label{fl}
 Let $f:\X \to S$ be the Faltings' example (Remark \ref{fal})
 and let $X$ be its generic fiber, which is a principally
polarized four-dimensional abelian variety over $\C(S)$. Since
$\X$ is not weakly isotrivial,
 it follows from Theorem \ref{iso} that the
 $\overline{\C(S)}/\C$-trace of
 $X$ is zero. Since $D_f \ne
 \End^0_S(\X)$, it follows that $\Is(X,\C(S))$ is {\sl infinite} and $\Is_1(X^2,\C(S),\ell)$ is
 also {\sl infinite} for all but finitely many primes $\ell$  congruent to $1$ modulo $4$.
 Recall that $G_s$ is connected (\ref{fal}) and therefore all
 endomorphisms of $X$ are defined over $\C(S)$, thanks to Remark
 \ref{endostab}. It follows from Corollary \ref{dim4ns} that $X$ is
 absolutely simple.
\end{ex}

We discuss isogeny classes of absolutely simple abelian fourfolds
in Section \ref{dim4}.

\section{Non-isotrivial Abelian schemes}
\label{AS}
\begin{proof}[Proof of Theorem \ref{iso}]
In one direction the assertion is almost obvious. Indeed, if
$\X_{S'}$ is a constant abelian scheme (resp. contains a non-zero
constant abelian subscheme) $W \times_ {\C} S'$ where $W$ is an
abelian variety over $\C$ then $X\times_{\C(S)} \C(S')$ is
isomorphic to $W\times_ {\C}\C(S')=(W\times_
{\C}\C(S))\times_{\C(S)}\C(S')$ (resp. contains an abelian
subvariety that is defined over $\C(S')$ and is isomorphic to
$W\times_ {\C}\C(S')$ over $\C(S')$. (By the way,  we did not use
an assumption that $S'\to S$ is \'etale.)

In the opposite direction, let $W$ be an abelian variety over
$\C$,  let $W_{\eta}$ be the constant abelian variety  $W\times_
{\C}\C(S)$ and  $\bar{u}: W_{\eta}\times_{\C(S)}\overline{\C(S)}
\hookrightarrow X\times_{\C(S)}\overline{\C(S)}$  an embedding of
abelian varieties over $\overline{\C(S)}$. Let $L=\C(S)(X_3)$ be
the field of definition of all points of order $3$ on $X$ and $S'$
is the normalization of $S$ in $L$. By Remark \ref{endostab}(i),
$S'\to S$ is a finite \'etale map of smooth irreducible curves,
$\C(S')=L$ and all points of order $3$ on $X$  are $L$-rational.
Clearly, all torsion points of $W_{\eta}$ are $\C(S)$-rational; in
particular, all points of order  $3$ on $W_{\eta}$ are defined
over $L$. It follows from results of \cite{Silverberg} that all
$\overline{\C(S)}$-homomorphisms between $W_{\eta}$ and $X$ are
defined over $\C(S')$; in particular, $\bar{u}$ is defined over
$\C(S')$, i.e.,  there exists an embedding
$u:W_{\eta}\times_{\C(S)}\C(S')\hookrightarrow
X\times_{\C(S)}\C(S')$ of abelian varieties over $\C(S')$ such
that $\bar{u}$ is obtained from $u$ by  extension of scalars from
$\C(S')$ to $\overline{\C(S)}$. Notice that the $\C(S')$-abelian
varieties $X\times_{\C(S)}\C(S')$ and
$W_{\eta}\times_{\C(S)}\C(S')=W\times_ {\C}\C(S')$ are generic
fibers of abelian $S'$-schemes $\X_{S'}$ and $W\times_{\C} S'$
respectively.
  It follows that $u$ extends
to a certain homomorphism of abelian $S'$-schemes $W\times_{\C}
S'\to \X_{S'}$, which we  denote by $u_{S'}$. If $u$  is an
isomorphism of generic fibers (the isotrivial case) then $u_{S'}$
is an isomorphism of the corresponding abelian schemes; in
particular, $\X_{S'}$ is a constant abelian scheme. If $u$ is not
an isomorphism, it is still a closed emdedding; in particular, the
image $Y$ of  $u:W\times_ {\C}\C(S')\hookrightarrow
X\times_{\C(S)}\C(S')$ is an abelian subvariety in
$X\times_{\C(S)}\C(S')$ and the abelian varieties $Y$ and
$W\times_ {\C}\C(S')$ are isomorphic over $\C(S')$. Let $\Y$ be
the schematic closure of the image $Y$ in $\X_{S'}$. It follows
from Corollary 6 on p. 175 of \cite{Neron} that $\Y$ is the
N\'eron model of $Y$ over $S'$. Since $W\times_ {\C}\C(S')\cong Y$
over $\C(S')$ and the N\'eron model of $W\times_ {\C}\C(S')$ over
$S'$ is $W\times S'$, we conclude that $\X_{S'}$ contains an
abelian subscheme isomorphic to the constant abelian scheme
$W\times_{\C}S'$.
\end{proof}

\section{Isogeny classes of abelians schemes}
\label{isab} Let $\X \to S$ be a polarized abelian scheme of positive relative
dimension $d$. Let us consider the category $Is(\X)$, whose objects are pairs
$(\Y,\alpha)$ that consist of an abelian scheme $g:\Y\to S$ and an isogeny
$\alpha:\Y \to \X$ of abelian schemes and the set of morphisms
$\Mor((\Y_1,\alpha_1), (\Y_2,\alpha_2)):=\Hom_S(\Y_1,\Y_2)$ for any pair of
objects $(\Y_1,\alpha_1)$ and $(\Y_2,\alpha_2)$. Let us consider the category
$Is_s(\X)$, whose objects are pairs $(\Lambda, i)$ that consist of a
$\pi_1(S,s)$-module
 $\Lambda$, whose additive group is isomorphic to $\Z^{2d}$ and an
 embedding $i:\Lambda\hookrightarrow H_1(\X_s,\Z)$ of
 $\pi_1(S,s)$-modules
 and the set of morphisms $\Mor_{Is_s(\X)}((\Lambda_1,i_1), (\Lambda_2,i_2))$ is the
set
 \[\{a \in \Hom_{\pi_1(S,s)}(\Lambda_1,\Lambda_2)\mid  \exists u\in
 \End^0(\X_s) \text{ such that } i_2 a= u i_1\}.\]

\begin{thm}
\label{fiber} The functor \[\Psi_s:Is(\X) \to Is_s(\X), \ (\Y,\alpha) \mapsto
(H_1(\Y_s,\Z),\ \alpha_s:H_1(\Y_s,\Z) \to H_1(\X_s,\Z)),\] \[(\gamma: \Y_1 \to
\Y_2) \mapsto \gamma_s\in
 \Hom({\Y_1}_s,{\Y_2}_s)\bigcap
\Hom_{\pi_1(S,s)}(H_1({\Y_1}_s,\Z),H_1({\Y_2}_s,\Z))\] is an equivalence of
categories.
\end{thm}

\begin{proof}
First, we need to check that $\gamma_s \in
\Mor_{Is_s(\X)}((H_1({\Y_1}_s,\Z),H_1({\Y_2}_s,\Z))$, i.e. there exists $u\in
\End^0(\X_s)$ such that ${\alpha_2}_s \gamma_s=u {\alpha_1}_s$. (Recall that
$\alpha_1: \Y_1 \to \X, \alpha_2:Y_2\to \X$ are isogenies of abelian schemes
and $\gamma_s \in \Hom({\Y_1}_s,{\Y_2}_s)$.) Clearly, both ${\alpha_1}_s:
{\Y_1}_s \to \X_s,\ {\alpha_2}_s:\Y_2\to \X_s$ are isogenies of abelian
varieties. Then $u:={\alpha_2}_s \gamma_s {\alpha_1}_s^{-1}\in \End^0(\X_s)$
satisfies ${\alpha_2}_s \gamma_s=u {\alpha_1}_s$. Second, the injectiveness and
surjectiveness of $\Hom_S(\Y_1,\Y_2)\to \Hom({\Y_1}_s,{\Y_2}_s)\bigcap
\Hom_{\pi_1(S,s)}(H_1({\Y_1}_s,\Z),H_1({\Y_2}_s,\Z))$ follow from the rigidity
lemma and Grothendieck's theorem (Sect. \ref{rigid}) respectively.

Let $(\Lambda,i_s)$ be an object of $Is_s(\X)$, i.e., a
$\pi_1(S,s)$-module
 $\Lambda$, whose additive group is isomorphic to $\Z^{2d}$ and
 an embedding $i_s:\Lambda\hookrightarrow H_1(\X_s,\Z)$ of
 $\pi_1(S,s)$-modules. In order to check the {\sl essential surjectiveness} of $\Psi_s$,
 we need to construct an abelian scheme
 $g:\Y\to S$, an isogeny $\alpha:\Y\to \X$ and an isomorphism
$\phi_s:\Lambda\cong H_1(\Y_s,\Z)$ of
 $\pi_1(S,s)$-modules such that $\phi_s=\alpha_s i_s$. In order to do
 that, recall that the $\pi_1(S,s)$-module $\Lambda$ defines a
 certain local system $\U$ of free $\Z$-modules of rank $2d$ on $S$,
 whose fiber at $s$
 coincides with  $\Lambda$.
 In addition, $i_s$ defines an embedding of local systems $i:\U
 \hookrightarrow
 R_1 f_*\Z$, whose fiber at $s$ coincides with our ``original" $i_s$. Rank arguments
 imply that the corresponding embedding \[i:\U\otimes {\Q} \to
 R_1 f_*\Z\otimes {\Q}=R_1f_*\Q\] is, in fact, an isomorphism. This
 allows us to provide $\U$ with the structure (induced by $R_1 f_*\Z$) of the holomorphic family of
 polarized Hodge structures of type $(-1,0)+(0,-1)$ \cite[Sect.
 4.4]{Deligne}. The equivalence of the category of polarized abelian
 schemes and the category of holomorphic families of
 polarized Hodge structures of type $(-1,0)+(0,-1)$ over $S$ \cite[Sect.
 4.4.2 and 4.4.3]{Deligne} (based on results of \cite{Borel}) implies that there exist an abelian
 scheme $g:\Y \to S$, a homomorphism $\alpha\in \Hom_S(\Y,\X)$ and an
 isomorphism of local systems $\phi_S:\U \cong R_1 g_*\Z$ such that
 $\alpha\psi_S=i$. Taking the fiber of the latter equality at $s$, we
 get the desired  $\alpha_s(\psi_S)_s=i_s$. Clearly,
 $\alpha_s\in \Hom(Y_s,X_s)$ induces an isomorphism
 $H_1(\Y_s,\Q)\cong H_1(\X_s,\Q)$ and therefore is an isogeny.
 Applying results of Section \ref{isogenySS},  we conclude that $\alpha$ is an isogeny of
 abelian schemes.
\end{proof}

\begin{rem}
The degree of $\alpha$ coincides with the index $[H_1(\X_s,\Z):i_s(\Lambda)]$.
\end{rem}

\section{Finite isogeny classes}
\label{fin}

 Let $L$ (resp. $L_{\Q}$) be the image of
the group algebra $\Z[\pi_1(S,s)]$ (resp. of $\Q[\pi_1(S,s)]$) in
$\End(H_1(\X_s,\Z))$ (resp. in $\End_{\Q}(H_1(\X_s,\Q))$) induced by the
monodromy representation. Clearly, $L$ is an order in the $\Q$-algebra
$L_{\Q}$.

It follows from Jackobson's density theorem and the semisimplicity of the
monodromy representation (over $\Q$) that $L_{\Q}$ is a semisimple $\Q$-algebra
that coincides with the centralizer $\End_{D_f}(H_1(\X_s,\Q))$ of $D_f$ and
\[\End_{L_{\Q}}(H_1(\X_s,\Q))=D_f.\] Clearly, each $\pi_1(S,s)$-stable
$\Z$-lattice in $H_1(\X_s,\Q)$ is an $L$-module, whose additive group is
isomorphic to $\Z^{2d}$.

\begin{sect}
 \label{ZS} It follows from the Jordan-Zassenhaus theorem
\cite[Theorem 26.4]{Reiner} that there are, up to an isomorphism,
only finitely many $L$-modules, whose additive group is isomorphic
to $\Z^{2d}$.
\end{sect}

\begin{proof}[Proof of Theorem \ref{finite}]
 Since $D_f=\End^0_S(\X)$, we have
\[D_f=\End^0_S(\X)\subset \End^0(\X_s)\subset \End_{\Q}(H_1(\X_s,\Q).\]

Thanks to Theorem \ref{fiber}, it suffices to check that the set of isomorphism
classes of objects in $Is_s(\X)$ is finite. Thanks to the Jordan-Zassenhaus
theorem (Sect. \ref{ZS}), it becomes an immediate corollary of the following
statement.

\begin{lem}
\label{categ} Suppose that $(\Lambda_1,i_1)$ and $(\Lambda_2,i_2)$
are objects in $Is_s(\X)$.

Then $i_1(\Lambda_1)$ and $i_2(\Lambda_2)$ are $L$-submodules in
$H_1(\X_s,\Z)$, whose additive group is isomorphic to $\Z^{2d}$. If the
$L$-modules $i_1(\Lambda_1)$ and $i_2(\Lambda_2)$ are isomorphic then
$(\Lambda_1,i_1)$ and $(\Lambda_2,i_2)$ are isomorphic.
\end{lem}

\begin{proof}[Proof of Lemma \ref{categ}]
Recall that $\Lambda_1$ and $\Lambda_2$ are $\pi_1(S,s)$-modules, whose
additive groups are isomorphic to $\Z^{2d}$ and $i_1:\Lambda_1 \hookrightarrow
H_1(\X_s,\Z), \ i_2:\Lambda_2 \hookrightarrow  H_1(X_s,\Z)$ are embeddings of
$\pi_1(S,s)$-modules. Clearly, $i_1:\Lambda_1\cong i_1(\Lambda_1), \
i_2:\Lambda_2\cong i_2(\Lambda_2)$ are isomorphisms of $\pi_1(S,s)$-modules.
Both $i_1(\Lambda_1)$ and $i_2(\Lambda_2)$ are  $L$-submodules in
$H_1(\X_s,\Z)$,

Assume that there exists an isomorphism of $L$-modules $:\alpha:
i_1(\Lambda_1)\cong i_2(\Lambda_2)$. Clearly, $\alpha$ is an isomorphism of
$\pi_1(S,s)$-modules. Extending $\alpha$ by $\Q$-linearity, we obtain an
isomorphism \[u:i_1(\Lambda_1)\otimes\Q\cong i_2(\Lambda_2)\otimes\Q\] of
$L\otimes\Q$-modules. Recall that $L\otimes\Q=L_{\Q}, \
i_1(\Lambda_1)\otimes\Q= i_2(\Lambda_2)\otimes\Q=H_1(\X_s,\Q)$. We have \[u \in
\End_{L_{\Q}}(H_1(\X_s,\Q))\bigcap \Aut_{\Q}(H_1(\X_s,\Q))= D_f\bigcap
\Aut_{\Q}(H_1(\X_s,\Q)).\] Since $D_f\subset\End^0(\X_s)$, we obtain that $u
\in \End^0(\X_s)$. Since $u\in\End^0(\X_s)\subset\End_{\Q}(H_1(\X_s,\Q))$ is an
automorphism of the finite-dimensional $\Q$-vector space $H_1(\X_s,\Q)$, we
have $u^{-1}\in \End^0(\X_s)$, i.e., $u \in\End^0(\X_s)^{*}$. Now if we put
\[a:=i_2^{-1}\alpha i_1: \Lambda_1 \to i_1(\Lambda_1)\to i_2(\Lambda_2)\to
\Lambda_2\] then $a \in \Hom_{\pi_1(S,s)}(\Lambda_1,\Lambda_2)$ is an
isomorphism of $\pi_1(S,s)$-modules and $i_2a=u i_1,\ i_1 a^{-1}=u^{-1} i_2$.
This implies that $a$ is an isomorphism of $(\Lambda_1,i_1)$ and
$(\Lambda_2,i_2)$.
\end{proof}
\end{proof}

\section{Infinite isogeny classes}
\label{inf}
 We deduce Theorem \ref{main} from the two following
auxiliary statements.

\begin{lem}
\label{linalg} Let $V$ be a finite-dimensional vector space over a
field $k$, let $\I_V$ be the identity automorphism of $V$, let
$A\subset \End_k(V)$ a $k$-subalgebra that contains $k\cdot\I_V$
and is isomorphic to a direct sum of matrix algebras over fields.
Let $B\subset A$ be a semisimple $k$-subalgebra that contains the
center of $A$ but  does not coincide with $A$.

Then there exists a proper subspace $W\subset V$ that enjoys the
following properties:
\begin{itemize}
\item[(i)] There does exist $a\in A$ such that $a(V)=W$.
\item[(ii)] There does not exist $b\in B$ such that $b(V)=W$.
\end{itemize}
\end{lem}

In order to state the next lemma, let us choose a prime $\ell$ and consider the
natural $\Q_{\ell}$-linear representation \[\pi_1(S,s)\to
\Aut(H_1(\X_s,\Q))\subset \Aut_{\Q_{\ell}}(H_1(\X_s,\Q)\otimes_{\Q}\Q_{\ell})
\subset\Aut_{\Q_{\ell}}(H_1(\X_s,\Q_{\ell}));\]
here $H_1(\X_s,\Q_{\ell})=H_1(\X_s,\Q)\otimes_{\Q}\Q_{\ell}$. Clearly, this
representation remains semisimple and the centralizer of $\pi_1(S,s)$ in
$\End_{\Q_{\ell}}(H_1(\X_s,\Q_{\ell}))$ coincides with
$D_f\otimes_{\Q}\Q_{\ell}$. Let us put $D=D_f$. We have
\[\End^0_S(\X)\otimes_{\Q}\Q_{\ell}\subset D\otimes_{\Q}\Q_{\ell}\subset
\End_{\Q_{\ell}}(H_1(\X_s,\Q_{\ell})).\] Clearly, $E\otimes_{\Q}\Q_{\ell}$ is
the center of $D\otimes_{\Q}\Q_{\ell}$ and lies in $C\otimes_{\Q}\Q_{\ell}$,
which is the center of $\End^0_S(\X)\otimes_{\Q}\Q_{\ell}$. It is also clear
that both $\Q_{\ell}$-algebras $D\otimes_{\Q}\Q_{\ell}$ and
$\End^0_S(\X)\otimes_{\Q}\Q_{\ell}$ are semisimple
 finite-dimensional.

\begin{lem}
\label{tate} Let $\ell$ be a prime, $W$ a proper $\pi_1(S,s)$-stable
$\Q_{\ell}$-vector space of $H_1(\X_s,\Q_{\ell})$ and
$d_0=\dim_{\Q_{\ell}}(W)$. Then one of the following two conditions holds:
\begin{enumerate}
\item There exists $u \in D\otimes_{\Q}\Q_{\ell}$ with
$u(H_1(\X_s,\Q_{\ell}))=W$. \item
 There exists a sequence $\{Z(m)\}_{m=1}^{\infty}$
of $\C(S)$-abelian varieties $Z(m)$ that enjoy the following
properties:
\begin{itemize}
\item[(i)] For each positive integer $m$ there exists a
$\C(S)$-isogeny $Z(m)\to X$ of degree $\ell^{(2d-d_0) m}$.
\item[(ii)] Let $Y$ be an abelian variety  over $\C(S)$ and
 let $M_Y$ be the set of positive integers $m$ such that $Z(m)$
 is isomorphic to $Y$ over $\C(S)$. Then $M_Y$ is either empty or
 finite.
 \end{itemize}
 \end{enumerate}
\end{lem}

\begin{rem}
The statement (and the proof) of Lemma \ref{tate} is inspired by
\cite[Prop. 1, pp. 136--137]{Tate}.
\end{rem}

\begin{proof}[Proof of Theorem \ref{main}(modulo Lemmas
\ref{linalg} and \ref{tate})] Let us apply Lemma \ref{linalg} to \[k=\Q_{\ell},
V=H_1(\X_s,\Q_{\ell}), A=D\otimes_{\Q}\Q_{\ell},
B=\End^0_S(\X)\otimes_{\Q}\Q_{\ell}.\] We conclude that there exist a proper
$\Q_{\ell}$-vector subspace $W$ in $H_1(\X_s,\Q_{\ell})$, an element $v\in
D\otimes_{\Q}\Q_{\ell}$ with $v(H_1(\X_s,\Q_{\ell}))=W$ but there does {\sl
not} exist an element $u\in \End^0_S(\X)\otimes_{\Q}\Q_{\ell}$ with
$u(H_1(\X_s,\Q_{\ell}))=W$. Since $v$ lies in $D\otimes_{\Q}\Q_{\ell}$, it
commutes with $\pi_1(S,s)$ and therefore $W$ is $\pi_1(S,s)$-stable. Now the
result follows from Lemma \ref{tate}.
\end{proof}

\begin{proof}[Proof of Lemma \ref{tate}]
Let us consider the $\pi_1(S,s)$-stable $\Z_{\ell}$-lattice
$H_1(\X_s,\Z_{\ell})=H_1(\X_s,\Z)\otimes\Z_{\ell}$ in the $\Q_{\ell}$-vector
space $H_1(\X_s,\Q_{\ell})$. Notice that for each positive integer $m$ we have
a canonical isomorphism of free $\Z/\ell^m\Z$-modules
\[H_1(\X_s,\Z_{\ell})/\ell^m H_1(\X_s,\Z_{\ell})=H_1(\X_s,\Z)/\ell^m
H_1(\X_s,\Z)=H_1(\X_s,\Z/\ell^m\Z)\] (induced by the canonical isomorphism
$Z/\ell^m\Z= \Z_{\ell}/\ell^m \Z_{\ell}$) that commutes with the actions of
$\pi_1(S,s)$.

The intersection $T:=W\bigcap H_1(\X_s,\Z_{\ell})$ is a $\pi_1(S,s)$-stable
free pure $\Z_{\ell}$-submodule of rank $d_0$
  in $H_1(\X_s,\Z_{\ell})$. Clearly,
$T\subset W$ and the natural map $T\otimes_{\Z_{\ell}}\Q_{\ell}\to W$ is an
isomorphism of $\Q_{\ell}$-vector spaces. The image $T_m$ of $T$ in
$H_1(\X_s,\Z/\ell^m\Z)$ is a free $Z/\ell^m\Z$-submodule of rank $d_0$. The
preimage of $T_m$ in $H_1(\X_s,\Z_{\ell})$ coincides with
$T+\ell^mH_1(\X_s,\Z_{\ell})$.

We write $\Lambda_m$ for the preimage of $T_m$ in $H_1(\X_s,\Z)$. Clearly,
$T+\ell^mH_1(\X_s,\Z_{\ell})$ contains $\Lambda_m$ and the natural map
$\Lambda_m\otimes\Z_{\ell} \to T+\ell^mH_1(\X_s,\Z_{\ell})$ is an isomorphism
of free $\Z_{\ell}$-modules. Notice that  $\Lambda_m$ is a $\pi_1(S,s)$-stable
subgroup of index $\ell^{m(2d-d_0)}$ in $H_1(\X_s,\Z)$ and contains $\ell^m
H_1(\X_s,\Z)$. It follows from Theorem \ref{fiber}  that there exists an
abelian scheme $h_m: \ZZ(m) \to S$  and an isogeny $\gamma(m): \ZZ(m)\to \X$ of
abelian schemes of degree $\ell^{m(2d-d_0)}$ such that
$\gamma(m)_s(H_1(\ZZ(m)_s,\Z))=\Lambda_m$. (Here $\ZZ(m)_s$ is the fiber of
$\ZZ(m)$ over $s$.) Since $\Lambda_m$ contains $\ell^m H_1(\X_s,\Z)$, there
exists an isogeny $\gamma(m)_s':\X_s \to \ZZ(m)_s$ of degree $\ell^{m d_0}$
such that the compositions
\[\gamma(m)_s\gamma(m)_s':\X_s \to \ZZ(m)_s \to
\X_s, \quad \gamma(m)_s'\gamma(m)_s: Z(m)_s \to \X_s\to \ZZ(m)_s\] coincide
with multiplication(s) by $\ell^m$. Clearly,
\[\gamma(m)_s(H_1(\ZZ(m)_s,\Z_{\ell}))=\gamma(m)_s(H_1(\ZZ(m)_s,\Z))\otimes\Z_{\ell}=
\Lambda_m\otimes\Z_{\ell}=T+\ell^m H_1(\X_s,\Z_{\ell}).\] We also have
$H_1(Z(m)_s,\Z)\supset\gamma(m)'_s(H_1(\X_s,\Z))\supset \ell^m
H_1(\ZZ(m)_s,\Z)$ and therefore
\[H_1(\ZZ(m)_s,\Z_{\ell})\supset\gamma(m)'_s(H_1(\X_s,\Z_{\ell}))\supset
\gamma(m)'_s \gamma(m)_s(H_1(\ZZ(m)_s,\Z_{\ell}))=\]
\[ \ell^m H_1(\ZZ(m)_s,\Z_{\ell}).\]
 Since $\gamma(m)_s \in \Hom_{\pi_1(S,s)} (H_1(\ZZ(m)_s,\Z),H_1(\X_s,\Z))$, we
have \[\gamma(m)_s^{-1} \in \Hom_{\pi_1(S,s)} (H_1(\X_s,\Q),H_1(Z(m)_s,\Q)).\]
It follows that $\gamma(m)_s'=\ell^m \gamma(m)_s^{-1}$ lies in
\[\Hom_{\pi_1(S,s)} (H_1(\X_s,\Q),H_1(\ZZ(m)_s,\Q))\bigcap \Hom
(H_1(\X_s,\Z),H_1(\ZZ(m)_s,\Z))=\]
\[\Hom_{\pi_1(S,s)}(H_1(\X_s,\Z),H_1(\ZZ(m)_s,\Z)).\] This implies that
$\gamma(m)_s'$ coincides with the ``fiber over" $s$ of a certain isogeny of
abelian schemes $\gamma(m)':\X \to \ZZ(m)$ of degree $\ell^{m d_0}$.  The
rigidity lemma implies that the compositions \[\gamma(m)\gamma(m)':\X \to
\ZZ(m) \to \X, \quad \gamma(m)'\gamma(m):\ZZ(m) \to \X\to \ZZ(m) \] coincide
with multiplication(s) by $\ell^m$.

 The generic fiber $Z(m)$ of $\ZZ(m)$ is a $\C(S)$-abelian variety and
$\gamma(m)$ and $\gamma(m)'$ induce  $\C(S)$-isogenies $\gamma(m)_{\eta}:Z(m)
\to X, \ \gamma(m)'_{\eta}:X \to Z(m)$. Their degrees are $\ell^{m(2d-d_0)}$
and $\ell^{m d_0}$ respectively. Their composition(s)
\[\gamma(m)_{\eta}\gamma(m)'_{\eta}: X \to Z(m) \to X, \
\gamma(m)'_{\eta}\gamma(m)_{\eta}:Z(m) \to X\to Z(m)\] coincide with
multiplication(s) by $\ell^m$.

Suppose that the condition (2) is not fulfilled. Then there exist an infinite
set $I$ of positive integers and  an abelian variety $Y$ over $\C(S)$ such that
$Y$ and $Z(i)$ are isomorphic over $\C(S)$ for all $i \in I$. Let $n$ be the
smallest element of $I$. For each $i\in I$ let us fix a $\C(S)$-isomorphism
$v_{i,\eta}:Z(n) \to Z(i)$. Clearly, $v_{i,\eta}$ extends to an isomorphism of
abelian schemes $\ZZ(n)\to\ZZ(i)$, which we  denote by $v_i$. Since $v_i$ is an
isomorphism, its specialization $v_{i,s}$ at $s$ satisfies
\[v_{i,s}(H_1(\ZZ(n)_s,\Z_{\ell}))=H_1(\ZZ(i)_s,\Z_{\ell}).\] Let us consider
the composition \[u_i:=\gamma(i)v_i\gamma(n)' : \X \to \ZZ(n)\to \ZZ(i) \to
\X.\] We have \[u_i \in \End_S(\X)\subset \End_S(\X)\otimes \Z_{\ell}\subset
\End_S(\X)\otimes \Q_{\ell} =\End^0_S(\X)\otimes_{\Q}\Q_{\ell}.\] Since
\[H_1(\ZZ(n)_s,\Z_{\ell})\supset \gamma(n)_s'(H_1(\X_s,\Z_{\ell}))\supset \ell^n
H_1(\ZZ(n)_s,\Z_{\ell}),\] we conclude that \[H_1(\ZZ(i)_s,\Z_{\ell})\supset
v_i\gamma(n)_s'(H_1(\ZZ(i)_s,\Z_{\ell}))\supset \ell^n
H_1(\ZZ(i)_s,\Z_{\ell})\] and therefore \[\ell^i
H_1(\X_s,\Z_{\ell})+T=\gamma(i)(H_1(\ZZ(i)_s,\Z_{\ell}))\supset
\gamma(i)v_i\gamma(n)'(H_1(\X_s,\Z_{\ell}))\supset \] \[\ell^n
\gamma(i)(H_1(\ZZ(i)_s,\Z_{\ell}))=\ell^n(\ell^i
H_1(\X_s,\Z_{\ell})+T)=\ell^{n+i}H_1(\X_s,\Z_{\ell})+\ell^n T.\] It follows
that for all $i\in I$, \[\ell^i H_1(\X_s,\Z_{\ell}))+T\supset
\gamma(i)v_i\gamma(n)'(H_1(\X_s,\Z_{\ell}))=u_i(H_1(\X_s,\Z_{\ell}))\supset\]
\[\ell^{n+i}H_1(\X_s,\Z_{\ell})+\ell^n T.\] Recall that
$\End_S(\X)=\End_{\C(S)}(X)$ is a free commutative
 group of finite rank and therefore the $\Z_{\ell}$-lattice
 $\End_S(\X)\otimes\Z_{\ell}$ in $\End_S(\X)\otimes_{\Q}\Q_{\ell}$
 is a compact metric space with respect to the $\ell$-adic topology. So, we can extract from $\{u_i\}_{i\in I}$ a
 subsequence $\{u_j\}_{j\in J}$ that converges to a limit
 \[u\in\End_S(\X)\otimes\Z_{\ell} \subset
 \End_{\Z_{\ell}}(H_1(\X_s,\Z_{\ell})).\]
 We may assume that there is a sequence of nonnegative integers $\{m_j\}_{j\in
 J}$ that tends to infinity and
 such that for all $j$,
 $u-u_j\in \ell^{m_j}\cdot \End_S(\X)\otimes\Z_{\ell}$.
 In particular,
$(u-u_j)(H_1(\X_s,\Z_{\ell}))\subset \ell^{m_j}\cdot H_1(\X_s,\Z_{\ell})$.
 It follows that
 \[u(H_1(\X_s,\Z_{\ell}))=\{\lim u_j(c)\mid c \in
 H_1(\X_s,\Z_{\ell})\}.\]
 This implies easily that $u(H_1(\X_s,\Z_{\ell}))\subset T$.
 On the other hand, if $t \in \ell^n T$ then for each $j\in J$
 there exists $c_j \in H_1(\X_s,\Z_{\ell})$ with $u_j(c_j)=t$.
 Since $H_1(\X_s,\Z_{\ell})$ is a compact metric space with respect to the $\ell$-adic topology,
 we can extract from $\{c_j\}_{j\in J}$ a
 subsequence $\{c_k\}_{k\in K}$ that converges to a limit $c \in
 H_1(\X_s,\Z_{\ell})$. Then
 $u(c)=\lim u_k(c) =\lim u_k(c_k)=t$.
 It follows that
 $\ell^n T \subset u(H_1(\X_s,\Z_{\ell}))\subset T$.
 This implies that
 $u(H_1(\X_s,\Q_{\ell}))=u(H_1(\X_s,\Z_{\ell}))\otimes_{\Z_{\ell}}\Q_{\ell}=
 T\otimes_{\Z_{\ell}}\Q_{\ell}=W$.
\end{proof}

\begin{proof}[Proof of Lemma \ref{linalg}]
{\bf Step 1}. {\sl Reduction to the case of simple} $A$. Suppose that the
semisimple $k$-algebra $A$ splits into a direct sum $A=A_1\oplus A_2$ of
non-zero semisimple $k$-algebras $A_1$ and $A_2$. Let $e_1$ and $e_2$ be the
identity elements of $A_1$ and $A_2$ respectively. Clearly, $\I_V=e_1+e_2, e_1
e_2=e_2e_1=0, e_1^2=e_1, e_2^2=e_2$. It is also clear that both $e_1$ and $e_2$
lie in the center of $A$ and therefore in the center of $B$. Let us put
\[V_1=e_1 V, V_2=e_2 V, \ B_1=e_1 B= B e_1\subset A_1, B_2=e_2 B= B e_2\subset
A_2;\] we have $V=V_1\oplus V_2, \  B = B_1 \oplus B_2$. In addition, $A_1$
acts trivially on $V_2$ and $A_2$ acts trivially on $V_1$. So, we may view
$A_1$ as a subalgebra of $\End_k(V_1)$ and $A_2$ as a subalgebra of
$\End_k(V_2)$ respectively. Obviously, the center of $A_i$ lies in the center
of $B_i$ for both $i=1,2$. Clearly, either $B_1 \ne A_1$ or $B_2 \ne A_2$. It
is also clear that the validity of the assertion of Lemma \ref{linalg} for
$(V_i,A_i,B_i)$ with $A_i\ne B_i$ implies its validity for $(V,A,B)$. It
follows that it suffices to prove Lemma \ref{linalg} under an additional
assumption that $A$ is simple. So, further we assume that $A$ is a simple
$k$-algebra and therefore is isomorphic to a matrix algebra $\M_r(E)$ of size
$r$ over a field $E$. Since $E$ is the center of $\M_r(E)$, it is a finite
algebraic extension of $k$. The field $E$ lies in the center of $D$, i.e., $D$
is a $E$-algebra. It follows easily that a semisimple $k$-algebra $D$ is a also
a semisimple $E$-algebra.

{\bf Step 2}. {\sl Reduction to the case of simple} $V=E^r$. Since all
$\M_r(E)$-modules of finite $k$-dimension are isomorphic to direct sums of
finite number of copies of the standard module $E^r$, we may assume that
$V=E^r, A=\M_r(E)=\End_E(V)$. Clearly, if $W$ is $k$-vector subspace in $V$
then one can find  $a \in \End(V)$ with $u(V)=W$ if and only if $W$ is a
$E$-vector subspace in $V$. So, in order to prove Lemma \ref{linalg}, it
suffices to prove the following statement.

\begin{prop}
\label{mat}
 Let $V$ be a vector space of finite dimension $r$ over a
field $E$. Let $B\subset \End_E(V)$ be a semisimple $E$-subalgebra
that contains $\I_V$. If $B \ne \End_E(V)$ then there exists a
proper $E$-vector subspace $W$ of $V$ that enjoys the following
property: there does not not exist $b \in B$ with $b(V)=V$.
\end{prop}

\begin{proof}[Proof of Proposition \ref{mat}]
Suppose that $B$ is not simple, i.e. $B$ splits into a direct sum $B=B_1\oplus
B_2$ of non-zero summands $B_1$ and $B_2$. Let $e_1$ and $e_2$ be the identity
elements of $A_1$ and $A_2$ respectively. Clearly $V$ splits into a direct sum
$V=e_1 V\oplus e_2 V$ of non-zero $B$-stable subspaces $e_1 V$ and $e_2 V$. It
is also clear that for each $b\in B$ the image \[bV=bV_1\oplus bV_2, \
bV_1\subset V_1, bV_2\subset V_2.\] In particular, if $W$ is {\sl not} a direct
sum of a subspace of
 $e_1 V$ and a subspace of $e_2 V$ then it does not coincide with
 $bV$ for any choice of $b \in B$. This proves Proposition
 \ref{mat}in the case of non-simple $B$. So further we assume that
 $B$ is a simple $E$-algebra. Let $F$ be the center of $B$.
 Clearly, $F$ is a field that contains $E\cdot\I_V$. If $F \ne E\cdot\I_V$
 then every $b(V)$ is an $F$-vector space; in particular, its
 $k$-dimension is divisible by the degree $[F: E\cdot\I_V]$, which is
 greater than $1$. So, none of $E$-subspaces $W$ of $E$-dimension
 $1$ is of the form $b(V)$. So, without loss of generality, we may
 assume that $F = E\cdot\I_V$, i.e., $B$ is a central simple
 $F$-algebra. Then there exist a division algebra $H$ of finite
 dimension over its center $E$ and a positive integer $m$ such
 that $B$ is isomorphic to the matrix algebra $\M_m(H)$. It follows from the
 classification of modules over central simple algebras that the
 left
 $B=\M_m(H)$-module $V$ is isomorphic to a direct sum of finitely
 many copies of the standard module $H^m$. Clearly, $V$ carries a
 natural structure of right $H$-module in such a way that every
 $b(V)$ is a right $H$-submodule. In particular, $\dim_E(bV)$ is
 divisible by $\dim_E(H)$. This implies that if $H \ne E$ then
 none of $E$-subspaces $W$ of $E$-dimension
 $1$ is of the form $bV$. So, without loss of generality, we may
 assume that $H=E$ and $B \cong M_m(E)$. Since $B \subset
 \End_E(V)$ but $B \ne\End_E(V)$, we conclude that $m<r$.
 Now dimension arguments
 imply that the $B=M_m(E)$-module $V=E^r$ is isomorphic to the
 direct sum of $r/m$ copies of the standard module $E^m$; in
 particular, $m\mid r$. This implies that $\dim_E(bV)$ is
 divisible by $r/m$. Since $m<r$, none of $E$-subspaces $W$ of $E$-dimension
 $1$ is of the form $bV$.
\end{proof}

\end{proof}

\section{Quaternions and SU(2)}
\label{append}

This stuff is (or should be) well-known. However, I was unable to
find a proper reference. (However, see \cite[Lecture 1, Example
12]{Postnikov}.)

Let $E$ be a field of characteristic zero, $F/E$ its quadratic extension. We
write $\sigma:z\mapsto \bar{z}$ for the only nontrivial $E$-linear automorphism
of $F$, which is an involution. Fix a non-zero element $\i\in E$ such that
$\sigma(\i)=-\i$. Clearly, \[F=E+E\cdot \i,\ 0\ne  -a:=\i^2\in E.\] Let $V$ be
a two-dimensional $E$-vector space. We write $\Aut_{F}^{1}(V)$ for the group of
$F$-linear automorphisms $u$ of $V$ with ${\det_F}(u)=1$. Here
${\det}_F:\Aut_F(V) \to F^{*}$ is the determinant map. If $\omega: V \times V
\to F$ is a map then we write $\Aut_{F}(V,\omega)$ for  the group of all
$F$-linear automorphisms $u$ of $V$  such that $\omega(ux,uy)=\omega(x,y) \
\forall \ x,y \in V$ and put
\[\Aut^{1}_{F}(V,\omega):=\Aut_{F}(V,\omega)\bigcap \Aut_{F}^{1}(V).\] Let
$\psi: V \times V \to F$ be a non-degenerate $F$-sesquilinear Hermitian form.
Let us reduce $H$ to a diagonal form, i.e., pick $e_1 \in V$ with
$\psi(e_1,e_1)\ne 0$ and let $V_2$ be the orthogonal complement in $V$ to $e_1$
with respect to $\psi$. The non-degeneracy of $\psi$ implies that $V_2$ is
one-dimensional and the restriction of $\psi$ to $V_2$ is also non-degenerate,
i.e., $\psi$ does not vanish on non-zero elements of $V_2$. It follows that if
$e_2$ is a non-zero element of $V_2$ then $\{e_1,e_2\}$ is an orthogonal basis
of $V$ then \[\psi(z_1 e_1+z_2 e_2, w_1 e_1+w_2 e_2)= b_1 z_1\overline{w_1}+b_2
z_2\overline{w_2} \ \forall \ z_1,z_2,w_1,w_2\in F\] where $0 \ne
b_1:=\psi(e_1,e_1)\in E, \  0 \ne b_2:=\psi(e_2,e_2)\in E$. Let us put
 $\discr(\psi)=b_1 b_2$.
It is known \cite[Ch. 9, \S 2]{Bourbaki} that the (multiplicative) class of
$\discr(\psi)$ in $E^{*}/N_{F/E}(F^{*})$ does not depend on the choice of
basis. Here $N_{F/E}: F \to E$ is the norm map. Let us put $b=b_2/b_1 \in
F^{*}$ and consider the Hermitian form \[\psi':=\frac{1}{b_1}\psi: V \times  V
\to F, \] \[\psi'(z_1 e_1+z_2 e_2, w_1 e_1+w_2 e_2)=  z_1\overline{w_1}+b
z_2\overline{w_2} \ \forall \ z_1,z_2,w_1,w_2\in F.\] Clearly, the classes of
$b$ and $\discr(\psi)$ in $E^{*}/N_{F/E}(F^{*})$ do coincide and
\[\Aut_{F}^{1}(V,\psi')=\Aut_{F}^{1}(V,\psi).\] Let us consider the {\sl cyclic}
algebra \cite[Sect. 15.1]{Pierce} \[D=(F,\sigma,-b)\cong
(F,\sigma,-\discr(\psi)).\] Recall that $D$ is the four-dimensional central
simple $E$-algebra that contains $F$ and coincides as (left) $F$-vector space
with $F\oplus F\cdot \j$ where $\j$ is an element of $A^{*}$ such that
$\j^2=-b, \ \j z \j ^{-1}=\sigma(z)=\bar{z}\ \forall z \in F$. If we put
$\k:=\i\j$ then \[A=F\oplus F\cdot \j=E\cdot 1\oplus E\cdot \i\oplus E\j\oplus
E\cdot \k\] with \[\i^2=-a,\j^2=-b, \k^2=-ab, \k=\i\j=-\j\i, \i\k=-a\j=-\k\i,
\k\j=-b\i=-\j\k.\] Clearly, $D=\left(\frac{-a,-b}{E}\right)$ (see \cite[Sect.
15.4]{Pierce}). Let us consider the standard $F$-linear involution $D \to D, \
q \mapsto \bar{q}$ that sends $\i,\j,\k$ to $-\i,-\j,-\k$ respectively.
Clearly,
\[\bar{q}q=q\bar{q}=x^2+ay^2+bs^2+abt^2\subset E \ \forall \ q=x\cdot
1+y\cdot\i+s\cdot\j+t\cdot\k; \ x,y,s,t\in E.\]
Let us consider the
``quaternionic" Hermitian $D$-sesquilinear form
\[\phi_D:D \times D \to D, (q_1,q_2) \mapsto q_1\overline{q_2}.\]
Taking the compositions of
$\phi_D:D\times D \to D$ with the projection maps
\[D=F\oplus F\cdot \j \twoheadrightarrow F, \ D=F\oplus F\cdot \j \twoheadrightarrow F \j,\]
we get
$E$-bilinear forms
\[H_D:D\times D \to F, \ A_D:D\times D \to F\]
defined by
\[\phi_D(q_1,q_2)=H_D(q_1,q_2)\cdot 1 +A_D(q_1,q_2)\cdot\j.\]
Clearly, both $H_D$
and $A_D$ are $F$-linear with respect to first argument. Since
\[\phi_D(q_1,q_2)=\overline{\phi_D(q_2,q_1)}, \ \phi_D(q,q)=q\bar{q}\in E\subset
F,\]
we have
\[H_D(q_1,q_2)=\overline{H_D(q_2,q_1)}, \ A_D(q,q)=0 \ \forall
q_1,q_2,q \in D.\]
This means that
 $H_D$ is an Hermitian
$F$-sesquilinear form on $D$ and $A_D$ is an alternating
$F$-bilinear form on $D$.

Clearly, $A_D$ is not identically zero and therefore is
non-degenerate, since the $F$-dimension of $D$ is $2$. This
implies that the {\sl symplectic} group $\Aut_F(D,A_D)$ coincides
with $\Aut_{F}^{1}(V)$
  and therefore
\[\Aut_F(D,\phi_D)=\Aut_F(D,H_D)\bigcap
\Aut_F(D,A_D)=\Aut_F(D,H_D)\bigcap\Aut_{F}^{1}(V)=\] \[\Aut_{F}^{1}(D,H_D).\]
Clearly, $F\cdot 1$ and $F\cdot\j$ are mutually orthogonal with respect to
$H_D$. On the other hand \[H_D(z_1\cdot 1,z_2\cdot 1)=z_1\overline{z_1},\
H_D(z_1\cdot \j,z_2\cdot \j)=z_1\j\overline{z_2\j}=z_1\j(-\j)\overline{z_2}= b
z_1 \overline{z_2}\ \forall \ z_1,z_2\in F.\] This implies that the isomorphism
of $F$-vector spaces \[\kappa:V\cong D, z e_1+w e_2 \mapsto z\cdot 1+ w \cdot
j\] is an isomorphism of {\sl Hermitian} $F$-vector spaces $(V,\psi')$ and
$(D,H_D)$, i.e., \[H_D(\kappa(v_1),\kappa(v_2))=\psi'(v_1,v_2)\ \forall
v_1,v_2\in V.\] For every $q\in D$ we denote by $R(q)$ the $F$-linear operator
$D \to D,  \ d \mapsto d\cdot q$. Clearly, \[R(1)=\I_{D},\
R(q_1q_2)=R(q_2)R(q_1), \ R(x q_1+y q_2)=x R(q_1)+y R(q_2)\] for all
$q_1,q_2\in D;\ x,y\in E$.
 I claim that
\[\Aut_F(D,\phi_D)=
\{R(q)\mid q\bar{q}=1\}.\] Indeed, if $q\bar{q}=1$ then
\[\phi_D(R(q)q_1,R(q)q_2)=\phi_D(q_1q,q_2q)=q_1q\overline{q_2q}=q_1q\bar{q}\
\overline{q_2} =q_1\overline{q_2}=\phi_D(q_1,q_2),\] i.e., $R(q)$ preserves
$\phi_D$.
 On the other hand, if a
$E$-linear automorphism $u$ of $D$ preserves $\phi_D$ then
$u(1)\overline{u(1)}=\phi_D(u(1),u(1))=\phi_D(1,1)=1\bar{1}=1\cdot
1=1$ and
 $u'=R(u(1))^{-1} u$ also preserves $\phi_D$ and
satisfies $u'(1)=1$. This implies that for all $q\in D$, \[q=q\cdot 1= q=q\cdot
\bar{1}= \phi_D(q,1)=\phi_D(u'(q),u'(1))= u'(q)\cdot \bar{1}=u'(q),\] i.e.,
$u'$ is the identity map $\I_D$ and therefore $u=R(u(1))$. Since
$\Aut_F(D,\phi_D)=\Aut_{F}^{1}(D,H_D)$, we conclude that
$\Aut_{F}^{1}(D,H_D)=\{R(q)\mid q\bar{q}=1\}$.

 Viewing $D$ as the
left $D$-module, we get an embedding $D \subset \End_{E}(D)$.
Clearly, (the algebra) $D$ (of left multiplications) commutes with
all {right multiplications} $R(q)$ and therefore with
$\Aut_{F}^{1}(D,H_D)$ in $\End_{E}(D)$.
\begin{lem}
\label{SUD} The centralizer ${\mathcal D}$ of
$\Aut_{F}^{1}(D,H_D)$ in $\End_{E}(D)$ coincides with $D$.
\end{lem}
\begin{proof}
Clearly, ${\mathcal D}$ contains $D$.  Let us pick {\sl non-zero} integers $n$
and $m$ such that $n^2+a\ne 0,\ m^2+b \ne 0$ and put \[q_1=\frac{(n^2-a)\cdot
1+ 2n\cdot \i}{n^2+a}, \ q_2=\frac{(m^2-b)\cdot 1+ 2m\cdot \j}{m^2+b}.\] One
may easily check that $q_1\overline{q_1}=1=q_2\overline{q_2}$. Clearly,
\[R(q_1)==\frac{n^2-a}{n^2+a}\I_V+\frac{2n}{n^2+a}R(\i),
R(q_2)==\frac{m^2-b}{m^2+b}\I_V+\frac{2m}{m^2+b}R(\j)\]
 and therefore ${\mathcal D}$ commutes with $R(q_1)$ and
$R(q_2)$. Since $n$ and $m$ do not vanish, we conclude that
${\mathcal D}$ commutes with right multiplications $R(\i)$ and
$R(\j)$. Clearly, $R(\j(R(\i)=R(\k)$. This implies that ${\mathcal
D}$ commutes with $R(\k)$ and therefore commutes with all $R(q)$
($q\in D$). Now if $u\in {\mathcal D}$ then let us put $z:=u(1)\in
D$. For all $q\in D$ we have
$u(q)=u(R(q)(1))=R(q)(u(1))=R(q)(z)=zq$, i.e., $u=z\in
D\subset\End_E(D)$.
\end{proof}

Recall that the Hermitian $F$-vector spaces $(V,\phi')$ and
$(D,H_D)$ are isomorphic. Taking into account that
$\Aut_{F}^{1}(V,\psi)=\Aut_{F}^{1}(V,\psi')$ and applying Lemma
\ref{SUD}, we obtain the following statement.

\begin{thm}
\label{cyclic} Let us view $V$ as a four-dimensional $E$-vector
space. Then the centralizer $D$ of $\Aut_{F}^{1}(V,\psi)$  in
$\End_E(V)$ is a central simple four-dimensional $E$-algebra
isomorphic to the four-dimensional central simple $E$-algebra
$D=\left(\frac{-a,-b}{E}\right)$.
\end{thm}

\section{Abelian fourfolds}
\label{dim4}

\begin{thm}
\label{dim4s} Let $S$ be a smooth irreducible algebraic curve over
$\C$ and $f:\X \to S$ a polarized abelian scheme  of relative
dimension $4$ with generic fiber $X$.

Suppose that $\X$ is not weakly isotrivial and $\Is(X,\C(S))$ is
infinite. Then:
\begin{itemize}

\item[(i)] $X$ is absolutely simple.

 \item[(ii)] The center $E$ of
$D_f$ is a real quadratic field.

\item[(iii)] $D_f$ is a quaternion division $E$-algebra that is
unramified at one infinite place of $E$ and ramified at the other
infinite place.

\item[(iv)] $\End_{\C(S)}^0(X)$ is either $E$ or a CM-field of
degree $4$.

\item[(v)]
$\End^0(X)$ is a CM-field of degree $4$.

\item[(vi)] Let $s\in S$  and assume that $G_s$ is connected. Then
$\End_{\C(S)}^0(X)=\End^0(X)$ is a CM-field of degree $4$.
\end{itemize}
\end{thm}

\begin{proof}
The absolute simplicity of $X$ follows from Corollary \ref{dim4ns}, which
proves (i). This implies that $\End_{\C(S)}^0(X)$ has no zero divisors. Since
$E$ is isomorphic to a subalgebra of $\End_{\C(S)}^0(X)$, we conclude that $E$
is a number field and $[E:\Q]$  divides $2\dim(X)=8$. Since $\Is(X,\C(S))$ is
infinite, it follows from Theorem \ref{main0} that $D_f \ne \End_S^0(\X)$. By
the last sentence of \ref{semisimple}, $D_f \ne E$, i.e., $D_f$ is a
non-commutative central simple $E$-algebra. Since $\dim_E(D_f)$ divides
\[\dim_{E}(H_1(\X_s,\Q))=\frac{\dim_{\Q}(H_1(\X_s,\Q))}{[E:\Q]}=\frac{8}{[E:\Q]},\]
we conclude that $8/[E:\Q]$ is {\sl not} square-free. It follows that $E$ is
either $\Q$ or a quadratic field. On the other hand, Deligne \cite[Prop.
4.4.11]{Deligne} proved that if $E$ is either $\Q$ or an imaginary quadratic
field then $D_f = \End_S^0(\X)$. It follows that $E$ is a real quadratic field
(which proves (ii)) and $\dim_E(D_f)=4$. It follows from \cite[Prop.
4.4.11]{Deligne} combined with the inequality $D_f \ne \End_S^0(\X)$ that $D_f$
is unramified at one infinite place of $E$ and ramified at the another one.
This rules out the possibility that $D_f$ is a matrix algebra over $E$. It
follows that $D_f$ is a quaternion division $E$-algebra, which proves (iii).

We have \[E \subset \End_S^0(\X)\subset D_f, \ \End_S^0(\X)\ne D_f.\] Since
$\dim_E(D_f)=4$, we conclude that either $\End_S^0(\X)=E$ or $\End_S^0(\X)$ is
a field of degree $4$ that contains $E$. Since there is an embedding $E
\hookrightarrow \R$ such that $D_f\otimes_E\R$ is the standard quaternion
$\R$-algebra, either $E=\End_S^0(\X)$ or $\End_S^0(\X)$ is a degree four field
that is {\sl not} totally real. Since $\End_S^0(\X)=\End_{\C(S)}^0(X)$, we
conclude that either $E=\End_S^0(\X)=\End_{\C(S)}^0(X)$ or
$\End_S^0(\X)=\End_{\C(S)}^0(X)$ is a CM-field of degree $4$, which proves
(iv).

Let $s\in S$ be a point in general position and assume that $G_s$ is connected.
Suppose that $\End_{\C(S)}^0(X)=E$. We need to arrive to a contradiction. It
follows from the first assertion of Corollary \ref{endoS} that
\[\End^0(\X_s)=\End_S^0(\X)=\End_{\C(S)}(X).\] This implies that $\End^0(\X_s)=E$
is a real quadratic field. It follows from \cite[4.2, p. 566]{MZ} that
$\Hdg(\X_s)$ is $\Q$-simple. It follows from the first assertion of  Corollary
\ref{endoS} that $\End_{\C(S)}(X)=\End_S^0(X)=\End_{\pi_1(S,s)}(H_1(\X_s,\Q))$.
Applying Theorem \ref{main0}, we conclude that $\Is(X,K)$ is finite. The
obtained contradiction proves that $\End_{\C(S)}^0(X)$ is a CM-field of degree
$4$. It follows from Remark \ref{endostab} that $\End_{\C(S)}^0(X)=\End^0(X)$.
This proves (vi).

We still have to prove that $\End^0(X)$ is a CM-field of degree $4$
 without assuming the connectedness of $G_s$. However, there exists a
 (connected) finite \'etale cover $S'\to S$ such that for $s'\in S'$ the group
 $G_{s'}$ attached to the abelian scheme $\X'=\X\times_{S}S'$ is connected
  and all endomorphisms of $X$ are defined over
 the field $L=\C(S')$ (see Subsect. \ref{cover}, Example \ref{min} and Remark \ref{endostab}).
 By Remark \ref{twist}, $\Is(X\times_{K}L,L)$ is also infinite. Applying
  already proven  (vi) to the generic fiber $X\times_{K}L$ of $\X'$, we conclude that
 $\End^0(X)=\End^0(X\times_{K}L)$ is a CM-field of degree $4$. This proves (v).
\end{proof}

\begin{ex}
Let $f:\X \to S$ be the Faltings' example (Sect. \ref{fal}) and
\ref{fl}) and let $X$ be its generic fiber. It follows from
Theorem \ref{dim4s} combined with arguments in Example \ref{fl}
that $\End_{\C(S)}^0(X)=\End^0(X)$ is a CM-field of degree $4$
that is a purely imaginary quadratic extension of a totally real
quadratic field $E$ and $D_f$ is a quaternion division $E$-algebra
that is unramified at one infinite place of $E$ and ramified at
the other infinite place. In fact, the quaternion algebra $D_f$
was the starting point of Faltings' construction \cite[Sect.
5]{Faltings0}.
\end{ex}

\begin{thm}
\label{foursuf} Let $S$ be a smooth irreducible algebraic curve
over $\C$ and $f:\X \to S$ a polarized abelian scheme  of relative
dimension $4$ with generic fiber $X$. Suppose that $\X$ is not
isotrivial and $\End^0(X)=\End_{\C(S)}^0(X)$. Assume, in addition,
that $F:=\End^0(X)$ is a CM-field of degree $4$. Let $\mu_F$ be
the (multiplicative) group of roots of unity in $F$ and $r_F$ its
order. Then:

\begin{itemize}
\item[(0)] $\X$ is not weakly isotrivial.
 \item[(i)] There exists a finite \'etale cover $S^0\to
S$ such that $\C(S')/\C(S)$ is a cyclic extension of degree
dividing $r_F$ and the set $\Is(X,\C(S^0))$ is infinite.
\item[(ii)] Let $s\in S$ and assume that $G_s$ is connected. Then
$\Is(X,\C(S))$ is infinite.
\end{itemize}
\end{thm}

\begin{proof}
Clearly, $X$ is an absolutely simple abelian variety. It follows
easily that $\X$ is not weakly isotrviial.

 Let us put $F:=\End_{\C(S)}^0(X)$. By assumption, it is
a purely imaginary quadratic extension of a real quadratic field $E$. Pick a
point $s$ in general position. By Corollary \ref{endoS}, $\End^0(\X_s)=F$. Let
us put $V:=H_1(\X_s,\Q)$; recall that $V$ is an $8$-dimensional $\Q$-vector
space. Clearly, $V$ carries a a natural structure of two-dimensional $F$-vector
space and $\Gamma_s \subset \Aut_F(V)$. We can do better, using the
polarization on $\X$ that induces a polarization on $\X_s$, whose Riemann form
gives rise to a non-degenerate alternating $\pi_1(S,s)$-invariant $\Q$-bilinear
form $\phi: V \times V \to \Q$. Since the {\sl complex conjugation} $e \mapsto
e'$ on $F$ is the only positive involution on the CM-field $F=\End^0(\X_s)$,
all Rosati involutions on $\End^0(\X_s)=F$ coincide with the complex
conjugation. This implies  that $\phi(ex,y)=\phi(x,e'y) \ \forall x,y\in V,
e\in F$ where $e \mapsto e'$ is the {\sl complex conjugation} on $F$. Pick a
non-zero element $\alpha \in E$ with $\alpha'=-\alpha$. Then there exists a
unique (non-degenerate) $F$-Hermitian form \[\psi: V \times V \to F\] such that
$\phi(x,y)=\tr_{F/\Q}(\alpha \psi(x,y)) \ \forall x,y\in V$ \cite[Sect.
9]{Deligne900}. Here $\tr_{F/\Q}: F \to \Q$ is the trace map.
 Since $\phi$ is $\pi_1(S,s)$-invariant, $\psi$ is also
$\pi_1(S,s)$-invariant and therefore is $\Gamma_s$-invariant. Let us consider
the unitary group \[\UU(V,\psi)\subset \GL(V)\] of the $F$-vector space $V$
relative to $\psi$. {\sl A priori} $U(V,\psi)$ is an algebraic group over $E$,
but we regard it as an algebraic $\Q$-group, i.e., take its Weil restriction
over $\Q$. In particular, $U(V,\psi)(\Q)=\Aut_{F}(V,\psi)$. Clearly, $\Gamma_s
\subset \Aut_{F}(V,\psi)=U(V,\psi)(\Q)$ and therefore $G_s \subset
\UU(V,\psi)$. The semisimplicity of $G_s^0$ implies that $G_s^0\subset
\SU(V,\psi)$ where $\SU(V,\psi)$ is the special unitary group  of the
$F$-vector space $V$ relative to $\psi$. As above, we view $\SU(V,\psi)$ as an
algebraic $\Q$-(sub)group (of $\UU(V,\psi)$). Clearly,
\[\SU(V,\psi)(\Q)=\Aut_{F}^{1}(V,\psi):=\{u\in \Aut_{F}(V,\psi)\mid
{\det}_F(u)=1\}.\]

Since $\dim_F(V)=2$, the centralizer $\End_{\SU(V,\psi)(\Q)}(V)$ of
$SU(V,\psi)(\Q)$ in $\End_{\Q}(V)$ (and even in $\End_E(V)$) is {\sl strictly
greater} than $F$ (see Theorem \ref{cyclic}). Since the group $\SU(V,\psi)(\Q)$
is dense in $\SU(V,\psi)$ with respect to Zariski topology \cite[Ch. 5, Sect.
18, Cor. 18.3]{BorelLin}, the centralizer $\End_{\SU(V,\psi)}(V)$ of
$\SU(V,\psi)$ in $\End_{\Q}(V)$  coincides with $\End_{\SU(V,\psi)(\Q)}(V)$ and
therefore is also {\sl strictly greater} than $F$. This implies that if
$G_s\subset \SU(V,\psi)$ then $D_f \ne \End^0(X)$. Applying Theorem
\ref{main0}, we obtain the following statement.

\begin{lem}
\label{su2} If $G_s\subset \SU(V,\psi)$ then $\Is(X,\C(S))$ is infinite.
\end{lem}
{\sl End of the proof of Theorem} \ref{foursuf}. Since $G_s^0\subset
\SU(V,\psi)$, the assertion (ii) follows readily. In order to prove the
assertion (i),  recall (Sect. \ref{deligneT}) that $\Gamma_s^0:=\Gamma_s\bigcap
G_s^0$ is a normal subgroup of finite index in $\Gamma_s$. Since $G_s^0\subset
\SU(V,\psi)$, it follows that
${\det}_F(\Gamma_s)\subset \mu_F$. This implies that if we put
$\Gamma^0:=\{g\in \Gamma_s\mid {\det}_F(g)=1\}$ then $\Gamma^0$ is a normal
subgroup in $\Gamma_s$ and the quotient $\Gamma_s/\Gamma^0$ is a finite cyclic
subgroup, whose order divides $r_F$. Using the construction of Section
\ref{cover} applied to $\Gamma^0$, we get a finite \'etale Galois cover $S^0
\to S$ with Galois group $\Gamma_s/\Gamma^0$ and abelian $S^0$-scheme
$\X^0=\X\times_S S^0$ such that if $s_0\in S^0$ lies above $s$ then the image
$\Gamma_{s_0}$ of the corresponding monodromy representation coincides with
$\Gamma^0$. (Here we identify $\X_s$ with $\X^0_{s_0}$ and $H_1(\X^0_{s_0},\Q)$
with $H_1(\X_{s},\Q)=V$.) Clearly, $\Gamma_{s_0}=\Gamma^0\subset \SU(V,\psi)$.
Applying Lemma \ref{su2} to
 the abelian $S^0$-scheme
$\X^0 =\X\times_S S^0$, we conclude that $\Is(X,\C(S^0))$ is infinite. In order
to finish the proof of (i), notice that the field extension $\C(S^0)/\C(S)$ is
normal and its Galois group coincides with $\Gamma_s/\Gamma^0$.
\end{proof}

\begin{rem} We keep the notation and assumptions of Theorem \ref{foursuf}.
\begin{enumerate}
\item
 $\Hdg(\X_s)=\UU(V,\psi)$ \cite[Sect. 7.5]{MZ}. It follows from Deligne's Theorem \ref{deligneH}
 that $G_s^0=\SU(V,\psi)$.
 \item  It follows from Theorem \ref{cyclic} that the centralizer $\End_{\SU(V,\psi)}(V)$
of $\SU(V,\psi)$ in $\End_{\Q}(V)$ is a four-dimensional central simple
 $E$-algebra $D'$ containing $F$. \item Suppose
that $G_s=G_s^0$. It follows that $D_f$ contains $D'$ and
therefore does not coincides with $F=\End_S^0(\X)$. However, the
center ${\mathcal E}$ of $D_f$ lies in $F$. Since $F\subset
D'\subset D_f$ and the center of $D'$ is $E$, we conclude that
${\mathcal E}\subset E$. This means that either ${\mathcal E}=\Q$
or ${\mathcal E}=E$, because $E$ is a quadratic field. In both
cases $D_f$ is a central simple ${\mathcal E}$-algebra. Applying
Proposition 4.4.11 of \cite{Deligne}, we conclude that ${\mathcal
E}\ne \Q$, i.e., ${\mathcal E}=E$. Since $\dim_{\Q}(D_f)$ must
divide $8$ and $\dim_{\Q}(D')=8$, we conclude that $D_f=D'$.
Applying  again the same Proposition, we conclude that $D_f$ is
ramified at one infinite place of $E$ and unramified at another
one. In particular, $D_f$ is not isomorphic to the matrix algebra
$\M_2(E)$, i.e., $D_f$ is a {\sl quaternion} (division)
 $E$-algebra.
\end{enumerate}
\end{rem}

\begin{ex}
\label{mainex} Recall that the $5$th cyclotomic field $F:=\Q(\mu_5)$ is a
CM-field of degree $4$. Notice that $\mu_F$ is a cyclic group of order $10$,
i.e., $r_F=10$. Let us consider a smooth genus zero affine curve
\[S=\A^1\setminus \{0,1\}=\P^1\setminus \{0,1,\infty\}\] with coordinate
$\lambda$ and a family $\CC\to S$  of genus four smooth projective curves over
$S$ defined by the equation \[y^5=x(x-1)(x-\lambda)\] and the corresponding
family of jacobians $f:\J\to S$. Clearly, $\J$ is a (principally) polarized
abelian $S$-scheme of relative dimension four, $\C(S)=\C(\lambda)$ is the field
of rational functions in one variable over $\C$ and the generic fiber $J$ of
$f$ is the jacobian of the $\C(\lambda)$-curve $y^5=x(x-1)(x-\lambda)$ of genus
four. It follows from \cite[Lemma 2.2 and Prop. 2.7]{JN} that $\J$ is {\sl not}
isotrivial and \[\End^0(J)=\End_{\C(\lambda)}^0(J)=F.\] This implies that $J$
is absolutely simple and therefore $\J$ is {\sl not} weakly isotrivial. By
Theorem \ref{foursuf}(i), there exists a cyclic extension $L/\C(\lambda)$ that
has degree $r$ dividing $10$, is unramified outside $\{0,1,\infty\}$ and such
that $\Is(J,L)$ is {\sl infinite}. Using Kummer theory, we obtain easily that
there are nonnegative integers $a$ and $b$ such that \[0\le a<r\le 10, \ 0\le
b<r\le 10\] and $L=\C(\lambda)(\sqrt[r]{\lambda^a (\lambda-1)^b})$. Clearly,
$\C(\lambda)\subset L \subset \C(\sqrt[10]{\lambda},\sqrt[10]{\lambda-1})$. It
follows from Remark \ref{twist} that
$\Is(J,\C(\sqrt[10]{\lambda},\sqrt[10]{\lambda-1}))$ is infinite. Now Corollary
\ref{principal} implies that
 $\Is_1(J^2,\C(\sqrt[10]{\lambda},\sqrt[10]{\lambda-1}),\ell)$ is {\sl infinite} for all but finitely
 many primes $\ell$  congruent to $1$ modulo $4$.

Notice that $\C(\sqrt[10]{\lambda},\sqrt[10]{\lambda-1})$ is the
field of the rational functions on the affine Fermat curve $S':
u^{10}-v^{10}=1, \ u \ne 0, v \ne 0$ with
$u=\sqrt[10]{\lambda},v=\sqrt[10]{\lambda-1}$ and $S' \to S, \
\lambda=u^{10}$ is the corresponding finite \'etale cover.
\end{ex}


\begin{thebibliography}{99}

\bibitem{BS} A. Borel, J-P.\ Serre,  Th\'eor\`emes de finitude
en cohomologie galoisienne, Comm. Math. Helv.  39 (1964) 111--164.

\bibitem{BorelLin} A. Borel, Linear algebraic groups (W.A. Benjamin, New York
Amsterdam, 1969).

\bibitem{Borel} A. Borel,   Some metric properties of arithmetic quotients of symmetric spaces
 and an extension theorem, J. Differential Geometry  6  (1972) 543--560.

 \bibitem{Neron} S. Bosch, W. L\"utkebohmert, M. Raynaud, N\'eron
 models (Springer-Verlag, Berlin Heidelberg New York, 1990).

 \bibitem{Bourbaki} N. Bourbaki, Alg\'ebre, Chap. IX (Hermann, Paris, 1959).

\bibitem{Deligne} P. Deligne,  Th\'eorie de Hodge.II, Publ.
Math. IHES  40 (1971) 5--58.

\bibitem{DeligneK3} P. Deligne, La conjecture de Weil pour
les surfaces K3, Invent. Math.  15 (1972), 206--226.

\bibitem{Deligne900} P. Deligne,  Hodge cycles on abelian
varieties (notes by J. S. Milne). Lecture Notes in Math. 900, Springer-Verlag,
Berlin, 1982, pp. 9--100.

\bibitem{Faltings0} G. Faltings,  Arakelov's theorem for abelian varieties, Invent. Math.
73 (1983) 337--347.

\bibitem{Faltings1} G. Faltings,  Endlichkeitss\"atze f\"ur
abelsche Variet\"aten \"uber Z\"ahlkorpern, Invent. Math.  73 (1983) 349--366;
Erratum  75 (1984) 381.

\bibitem{Faltings2} G. Faltings,  `Complements to Mordell',
Chapter VI in: G. Faltings, G. Wustholz et al., Rational points. Third edition.
Aspects of Mathematics,  E6. Friedr. Vieweg \& Sohn, Braunschweig, 1992.

\bibitem{Grothendieck} A. Grothendieck,  Un th\'eor\`eme sur
les homomorphismes de sch\'emas ab\'elienes, Invent. Math.   2 (1966) 59--78.

\bibitem{JN} J. de Jong and R. Noot, `Jacobians with complex
multiplications',  Arithmetic algebraic geometry (eds G. van der Geer, F. Oort
and J. Steenbrink), Progress in Math., vol.  89 (Birkh\"auser, 1991), pp.
177--192.

\bibitem{Katz} N. Katz,  Nilpotent connection and monodromy
theorem: Application of a result of Turrutin, Publ. Math. IHES  39 (1971)
175--232.

\bibitem{Lang} S. Lang, Abelian varieties. Second edition
(Springer-Verlag, New York, 1983).

\bibitem{LOZ} H. W. Lenstra, F. Oort, Yu. G. Zarhin,  Abelian
subvarieties, J. Algebra  180 (1996) 513--516.

\bibitem{Masser} D. Masser,  Specializations of endomorphism
rings of abelian varieties, Bull. Soc. math. France  124 (1996) 457--476.


\bibitem{MZ} B.J.J.~Moonen, Yu.G.~Zarhin,  Hodge and Tate classes on simple
abelian fourfolds, Duke Math. J. 77 (1995) 553--581.

\bibitem{MZ2} B.J.J.~ Moonen, Yu. G.~Zarhin,  Hodge classes on abelian varieties of low
dimension, Math. Ann.  315 (1999) 711--733.

\bibitem{MB} L. Moret-Bailly, Pinceaux de vari\'et\'es ab\'eliennes, Ast\'erisque,
vol. 129 (1985).

\bibitem{MumfordH} D. Mumford, A note of Shimura's paper
``Discontinuous groups and abelian varieties", Math. Ann.  181 (1969) 345--351.

\bibitem{GIT} D. Mumford, J. Fogarty, F. Kirwan, Geometric
invariant theory, Third enlarged edition  (Springer-Verlag, Berlin Heidelberg
New York, 1994).

\bibitem{Mumford} D. Mumford, Abelian varieties, Second edition,
 Oxford University Press, London, 1974.

 \bibitem{Parshin} A. N.~Parshin, Minimal models of curves of genus 2
 and homomorphisms of elliptic curves over function fields of
 finite characteristic, Izv. AN SSSR ser. matem.  36 (1972) 67--109; Math.
 USSR Izv.  6 (1972) 60--108.

 \bibitem{Pierce} R. Pierce, Associative Algebras
 (Springer-Verlag, Berlin Heidelberg New York, 1982).

 \bibitem{Postnikov} M.M.Postnikov, Lie groups and Lie algebras  (MIR, Moscow,
 1986).

 \bibitem{Reiner} I. Reiner,  Maximal orders  (Academic Press,
 London, 1966).


\bibitem{Serre} J-P.\ Serre, Abelian $\ell$-adic representations
 and elliptic curves. Second edition (Addison-Wesley, New York, 1989).

 \bibitem{SerreTate} J-P.\ Serre, J.\ Tate,
 Good reduction of abelian varieties, Ann.\ of Math.\ (2)
 88  (1968) 492--517.

 \bibitem{Silverberg} A. Silverberg, Fields of definition for homomorphisms of abelian
varieties, J. Pure Applied Algebra  77 (1992) 253--262.

\bibitem{SZcomp} A. Silverberg,  Yu. G.~Zarhin,  Connectedness
results for $\ell$-adic representations associated to abelian varieties, Comp.
Math.  97 (1995) 273--284.

\bibitem{Tate} J. Tate, Endomorphisms of Abelian varieties over finite
fields, Invent. Math.  2 (1966) 134--144.

\bibitem{ZarhinFA} Yu. G.~Zarhin,  A finiteness theorem for isogenies of Abelian varieties over
function  fields of finite characteristic, Funktsional Anal. i Prilozh.  8
(1974),  No. 4, 31--34; Functional Anal. Appl.  8 (1974) 301--303.

\bibitem{ZarhinIz} Yu. G.~Zarhin,  Endomorphisms of Abelian varieties over fields of finite
characteristic,. Izv. Akad. Nauk SSSR ser. matem.  39 (1975) 272--277; Math.
USSR Izv.  9 (1975) 255--260.

\bibitem{ZarhinMZ1} Yu. G.~Zarhin,  Abelian varieties in characteristic
$P$, Mat. Zametki  19 (1976) 393--400; Mathematical Notes  19 (1976) 240--244.

\bibitem{ZarhinMZ2} Yu. G.~Zarhin,   Endomorphisms of Abelian varieties and points of finite order
in characteristic $P$, Mat. Zametki  21 (1977) 737--744; Mathematical Notes
 21 (1978) 415--419.

\bibitem{ZarhinIzv} Yu. G.~Zarhin, Weights of simple Lie algebras in the cohomology of
algebraic varieties, Izv. Akad. Nauk SSSR Ser. Mat. 48 (1984) 264--304; Math.
USSR Izv. 24 (1985) 245 - 281.

\bibitem{ZarhinIn} Yu. G.~Zarhin,  A finiteness theorem for unpolarized Abelian varieties over
number fields with prescribed places of bad reduction, Invent. Math.  79 (1985)
309--321.

\bibitem{ZarhinP} Yu. G.~Zarhin, A. N.~Parshin,  Finiteness problems in Diophantine
geometry,  Amer. Math. Soc. Transl. (2) 143 (1989) 35--102.
\end{thebibliography}
\end{document}